\documentclass{amsart}

\usepackage{amsmath,amssymb,amsthm,amscd,amsxtra}
\usepackage{amsfonts,amsbsy,mathrsfs}
\usepackage{mathtools} 
\usepackage{graphicx,subfigure}
\usepackage{epstopdf}
\usepackage{float}
\usepackage{multirow}
\usepackage{tikz}
\usetikzlibrary{positioning}
\usepackage{xcolor}
\usepackage{marginnote}
\usepackage{sidenotes}
\usepackage{changes}
\usepackage{aliascnt}
\usepackage{hyperref}
\usepackage{cleveref} 

\allowdisplaybreaks[3]

\newcommand{\bib}[1]{\bibitem{#1}}
\newcommand{\ignore}[1]{}

\theoremstyle{plain}

\newtheorem{theorem}{Theorem}[section]

\newaliascnt{lemma}{theorem}
\newtheorem{lemma}[lemma]{Lemma}
\aliascntresetthe{lemma}

\newaliascnt{proposition}{theorem}

\aliascntresetthe{proposition}

\newaliascnt{corollary}{theorem}
\newtheorem{corollary}[corollary]{Corollary}
\aliascntresetthe{corollary}

\newaliascnt{definition}{theorem}
\newtheorem{definition}[definition]{Definition}
\aliascntresetthe{definition}

\newaliascnt{assumption}{theorem}
\newtheorem{assumption}[assumption]{Assumption}
\aliascntresetthe{assumption}

\theoremstyle{definition}

\newaliascnt{example}{theorem}
\newtheorem{example}[example]{Example}
\aliascntresetthe{example}

\newaliascnt{algorithm}{theorem}

\aliascntresetthe{algorithm}

\newaliascnt{exercise}{theorem}

\aliascntresetthe{exercise}

\theoremstyle{remark}

\newaliascnt{remark}{theorem}
\newtheorem{remark}[remark]{Remark}
\aliascntresetthe{remark}

\numberwithin{equation}{section}

\Crefname{figure}{Figure}{Figures}

\DeclareMathOperator{\Real}{Re}

\DeclareMathOperator{\supp}{supp}

\DeclarePairedDelimiterX{\duality}[2]{\langle}{\rangle}{#1,#2}

\newcommand{\ii}{\mathrm{i}}
\newcommand{\II}{\mathcal I}
\newcommand{\RR}{\mathcal R}

\newcommand{\PP}{\mathcal P}
\newcommand{\BB}{\mathcal B}

\newcommand{\dd}{\,\mathrm{d}}
\newcommand{\xx}{\mathbf{x}}

\newcommand{\TT}{\mathcal{T}}
\newcommand{\faces}{\triangle}
\newcommand{\ifaces}{\mathring \triangle}

\newcommand{\A}{\mathcal A^\gamma}
\newcommand{\aA}{\mathcal A^{\gamma,\star}}
\newcommand{\At}{\mathcal A_T^\gamma}

\newcommand{\IdC}{\mathcal C_\infty^\gamma}
\newcommand{\aIdC}{\mathcal C_\infty^{\gamma,\star}}
\newcommand{\IdtC}{\mathcal C_{\infty,T}^\gamma}
\newcommand{\aIdtC}{\mathcal C_{\infty,T}^{\gamma,\star}}
\newcommand{\PtC}{\mathcal C_{L,T}^\gamma}
\newcommand{\aPtC}{\mathcal C_{L,T}^{\gamma,\star}}
\newcommand{\PC}{\mathcal C_L^\gamma}
\newcommand{\aPC}{\mathcal C_L^{\gamma,\star}}

\newcommand{\green}{\mathcal G_\infty^\gamma}
\newcommand{\agreen}{\mathcal G_\infty^{\gamma,\star}}

\newcommand{\msI}[1]{#1_{H,{\rm ms}}^\gamma}
\newcommand{\pmsI}[1]{#1_{H,L,{\rm ms}}^\gamma}
\newcommand{\apmsI}[1]{#1_{H,L,{\rm ms}}^{\gamma,\star}}

\newcommand{\abs}[1]{\left|#1\right|}
\newcommand{\jump}[1]{\left[#1\right]}
\newcommand{\norm}[1]{\left\lVert#1\right\rVert}
\newcommand{\tnorm}[1]{%
  \left\lvert\!\left\lvert\!\left\lvert #1
\right\rvert\!\right\rvert\!\right\rvert%
}
\newcommand{\gsnorm}[2]{\abs{#1}_{\gamma,#2}}
\newcommand{\ogsnorm}[1]{\abs{#1}_{\gamma}}
\newcommand{\gnorm}[2]{\norm{#1}_{\gamma,#2}}
\newcommand{\ognorm}[1]{\norm{#1}_{\gamma}}

\newcommand{\patch}{N^L(T)}
\newcommand{\pker}{W_h(\patch)}

\newcommand{\msII}[1]{#1^{\rm ms}_\infty}
\newcommand{\amsII}[1]{#1^{{\rm ms},\star}_\infty}
\newcommand{\pmsII}[1]{#1^{\rm ms}_L}
\newcommand{\apmsII}[1]{#1^{{\rm ms},\star}_L}
\newcommand{\Q}{\pmsII{\mathcal Q}}
\newcommand{\aQ}{\apmsII{\mathcal Q}}

\newcommand{\pmsIII}[2]{#1_{H,L,{\rm ms}}^{#2}}
\newcommand{\apmsIII}[2]{#1_{H,L,{\rm ms}}^{#2,\star}}

\begin{document}

\title[CIP-PML-LOD for Helmholtz problems]{A Continuous Interior Penalty Localized Orthogonal Decomposition Method for Large-Wavenumber Helmholtz Problems with UPML}

\author{Huangxin Chen}
\address{School of Mathematical Sciences and Fujian Key Laboratory on Mathematical Modeling and High Performance Scientific Computing, Xiamen University, Fujian, 361005, China}
\email{chx@xmu.edu.cn}

\author{Zihao Su}
\address{Department of Mathematics, Nanjing University, Jiangsu, 210093, China}
\email{652024210011@smail.nju.edu.cn}

\author{Dinghang Tan}
\address{School of Mathematical Sciences and Fujian Key Laboratory on Mathematical Modeling and High Performance Scientific Computing, Xiamen University, Fujian, 361005, China}
\email{dhtan@stu.xmu.edu.cn}

\author{Haijun Wu}
\address{Department of Mathematics, Nanjing University, Jiangsu, 210093, China}
\email{hjw@nju.edu.cn}



\begin{abstract}
We develop a continuous interior penalty localized orthogonal decomposition (CIP-LOD) method for large-wavenumber Helmholtz scattering problems. The exterior problem is truncated by a Cartesian PML, whereas CIP-stabilization is incorporated into the fine-scale reference discretization to mitigate pollution effect. At the same time, we then construct a transfer operator to unify or simplify the proofs of stability results and error estimates in the considered general variational multiscale methods. Under the resolution condition $\kappa H \leq C$, we derive the error estimates for the CIP-LOD solution $u_{h}$ and CIP-LOD solution $u_{H,L,ms}^{\gamma}$ explicitly, namely \[\norm{u_h - \pmsI{u}}_\gamma \leq O(e^{-\delta L/2\tau}) \norm{f} + O(H^{j+1})\abs{f}_{H^j(\Omega)}, 0 \le j \le 2.\] We have already developed an improved CIP-LOD method, in which the CIP penalty term is omitted during the computation of the local subproblems. In particular, we use the transfer operator to establish the stability and error estimates for that improved CIP-LOD method. Both theory and numerical results show that the improved CIP-LOD algorithm saves computational costs for local problems while keeping the relative error of the CIP-LOD solution clean. Numerical experiments confirmed our theoretical results, and when the right penalty parameter is chosen, the pollution error can be greatly reduced. With a purely imaginary penalty added to the CIP, CIP-FEM suppresses the pollution error less effectively than with the optimal penalty parameter, but still achieves a noticeable reduction compared with standard FEM. An interesting experimental observation is that the LOD resolution condition can be substantially relaxed in this setting, which greatly reduces the number of subproblems. 
\end{abstract}

\maketitle

\section{Introduction}

The numerical simulation of time-harmonic waves at high frequencies is a central challenge in acoustics, electromagnetics, geophysics, and scattering theory. With the time convention $e^{-\mathrm{i}\omega t}$, the frequency-domain model contains the indefinite Helmholtz operator $-\Delta-\kappa^2$. Standard conforming finite element methods (FEMs) face two main limitations: the mesh must resolve the wavelength, and the discrete phase error leads to the well-known pollution effect. Classical analyses for the $h$- and $hp$-FEM indicate that, in the preasymptotic regime, the error contains wave-number dependent terms beyond the best approximation error. This pollution effect is already visible in dispersion studies and remains a central issue in wave-number explicit finite element theory; see, for instance, \cite{Babuska1997,Deraemaeker1999,Esterhazy2012,Ihlenburg1995,Ihlenburg1997,Melenk2010,Melenk2011}. More recent analyses give sharper relative-error, $hp$-FEM, variable-coefficient, and piecewise-smooth-media results \cite{Bernkopf2024,ChaumontFreletNicaise2020,Galkowski2025,Graham2020,Lafontaine2022Sharp,Spence2023}. The difficulty is especially relevant for exterior radiation problems, where the oscillatory field must be approximated together with a reliable numerical realization of the condition at infinity.

Classical absorbing and radiation boundary conditions provide local approximations of the Sommerfeld condition \cite{Bayliss1980,Engquist1977,Hagstrom1999}. Perfectly matched layers (PMLs) was originally introduced for electromagnetic waves and subsequently interpreted through complex coordinate stretching. It incorporates a bounded artificial absorbing layer to absorb outgoing wave \cite{Berenger1994,Chew1994,Collino1998,Lassas2001}. Finite PML approximations and Cartesian uniaxial PMLs have been analyzed for acoustic, electromagnetic, and layered-media scattering problems \cite{Bramble2007,Bramble2013,Chen2010}. More recent high-frequency analyses show that PML truncation is exponentially accurate under suitable assumptions and is compatible with wave-number explicit finite element theory \cite{ChaumontFreletGallistlNicaiseTomezyk2022,Galkowski2023,Galkowski2024}. In particular, $hp$-FEM with PML truncation can be made pollution-free under appropriate resolution conditions \cite{Galkowski2024}. Continuous interior penalty finite element methods (CIP-FEMs) add penalties on jumps of normal derivatives across interelement faces. For Helmholtz problems with large wave number, preasymptotic estimates for linear, higher-order, and $hp$ CIP-FEMs clarify how penalty terms influence stability and pollution \cite{Burman2016,Du2015,Wu2014,Zhu2013}. Dispersion analysis and optimized penalty parameters further show that the pollution contribution can be reduced in practical regimes \cite{Li2024,Zhou2023}. For radiation problems, CIP-FEM can be combined with truncated Dirichlet-to-Neumann maps or PML truncation \cite{Li2019,Li2025,Li2020}. Error estimates have been established for both approaches. Related discontinuous Galerkin formulations provide another stabilized route for highly indefinite Helmholtz problems \cite{MelenkParsaniaSauter2013}. Localized orthogonal decomposition (LOD) was developed as a numerical homogenization framework for elliptic multiscale problems. The fine-scale information is compressed into localized coarse basis corrections \cite{Engwer2019,Henning2013,Malqvist2014}; see also the monograph \cite{Malqvist2020}. For high-frequency Helmholtz equation, LOD and multiscale Petrov--Galerkin methods build fine-scale trial or test correctors to encode subscale phase information. Stable and quasi-optimal behavior can be obtained under a natural coarse resolution condition of the form $\kappa H \lesssim 1$ in homogeneous settings \cite{Gallistl2015,Peterseim2017}. Extensions to heterogeneous or two-scale Helmholtz models have also been developed \cite{Brown2017,Ohlberger2017}. Multi-resolution and super-localized variants reduce or reorganize the corrector computations \cite{Freese2024,Hauck2022}, while multiscale generalized finite element methods provide related wave-number explicit results for heterogeneous media \cite{MaAlberScheichl2023}. These three strands are complementary. PML addresses the continuous truncation of the unbounded domain, CIP controls fine-grid dispersion and jump contributions, and LOD compresses fine-scale wave responses into localized coarse spaces. Therefore, this paper proposes a CIP-LOD method, i.e., an LOD procedure specifically designed for the CIP-FEM method. In this way, we can combine the advantages of CIP-FEM, UPML truncation, and the LOD method to efficiently handle large-scale problems while achieving good control of the pollution error.

In this work, we systematically analyze the influence of the UMPL parameters on the stability and error estimates for the localization and CIP-LOD. Although \cite{altmann2021numerical}\cite{Peterseim2017} provides a proof framework for the Helmholtz LOD problem with impedance boundary conditions, the influence of the UPML parameters has not been systematically treated. This influence can be partially understood as a high-contrast effect (see \cite{Babuska1997}). The proofs for the LOD method employed in this paper are all developed in the fully discrete setting, and thus can be viewed as a complement to \cite[Section~6]{Peterseim2017}. That section provides a general proof framework in the fully discrete setting, while this section elaborates on the relevant technical details. We also show that higher regularity of the right hand side $f$ improves the convergence rate of the idealized LOD error, it can be regarded as a supplement to \cite[Theorem~4.5]{Peterseim2017}. We then establish a pollution-free error estimate between the CIP-LOD and CIP-FEM solutions with the same penalty parameter. This makes CIP-LOD an efficient approach to solving large-scale CIP-FEM problems while retaining its pollution-suppression benefits. Moreover, we show that a purely imaginary penalty improves the stability of the LOD subproblems while preserving exponential decay. This relaxes the resolution condition on $\kappa H$. We rigorously prove this result for the Helmholtz problem with impedance boundary conditions. Numerical experiments demonstrate a similar effect for UPML problems.

Furthermore, we establish a unified abstract framework and use it to improve and extend the method. In the framework, a stability-transfer operator is introduced to improve the stability and error estimates. It applies to various variational multiscale methods, including constraint energy minimizing generalized multiscale FEM (CEM-GMsFEM), LOD, and Super-LOD. Using this framework, we improve CIP-LOD by omitting the CIP penalty in the local subproblems, which lowers the computational cost. We also prove the relative error between the improved method and the original CIP-LOD method. Moreover, the abstract framework applies to continuous, discontinuous-conforming, discontinuous-continuous, $\mathbf H(\rm curl)$, and $\mathbf H(\rm div)$ problems.

The remainder of this paper is organized as follows. \Cref{sec:model} introduces the PML-truncated model problem and the notation for the associated variational formulation. \Cref{sec:cip} presents the fine-scale CIP-PML discretization and introduces the quasi-interpolation operator. In \Cref{sec:lod}, we construct the ideal and localized CIP-LOD methods and prove exponential decay and localization estimates for the corrections. \Cref{sec:framework} develops an abstract stability-transfer framework. This framework is applied in \Cref{sec:alg-opt} to the improved method, in which the CIP penalty is omitted from the local corrector problems. Numerical experiments in \Cref{sec:numeric} illustrate the convergence, pollution-reduction, and localization behavior of the proposed method. Finally, \Cref{sec:con} concludes the paper.

\section{Model setting}\label{sec:model}
Let $\kappa>0$ denote the wave number and use the time factor $e^{- \ii \omega t}$. The paper starts from an exterior Helmholtz radiation problem in the whole space. The physical region is the square $B_{R_0}:=(-R_0,R_0)^2\subset\mathbb{R}^2$. The source $f^{\rm phys}$ is assumed to be compactly supported in $B_{R_0}$. The outgoing free-space solution $u^{\rm ext}$ is defined by
\begin{equation*}
	-\Delta u^{\rm ext}-\kappa^2 u^{\rm ext}=f^{\rm phys},\qquad\text{in }\mathbb{R}^2,
\end{equation*}
with the two-dimensional Sommerfeld radiation condition
\begin{equation*}
	\lim_{r\to\infty}r^{1/2}\big(\partial_r u^{\rm ext}-\ii \kappa u^{\rm ext}\big)=0,
	\qquad r=|x|.
\end{equation*}
The Sommerfeld condition is imposed at infinity and is not directly suitable for a finite element discretization. We therefore replace the exterior problem by a bounded-domain PML problem, which is the continuous model discretized in the subsequent sections.

Choose $R>R_0$. Set $B_R:=(-R,R)^2$, $B_{\rm pml}:=B_R\setminus\overline{B_{R_0}}$, and $d_{\rm pml}:=R-R_0$. The outer boundary is denoted by $\Gamma_{\rm out}:=\partial B_R$, and a homogeneous Dirichlet condition is imposed on $\Gamma_{\rm out}$. In the localization arguments below we also use the shorthand $\Omega:=B_R$. For $j=1,2$, let $\sigma_j=\sigma_j(x_j)$ be a nonnegative absorption profile that vanishes in the physical region and is bounded in the layer. A typical choice is the polynomial profile
\begin{align*}
	\sigma_j(x_j)=
	\begin{cases}
		0, & |x_j|\le R_0,\\[1mm]
		\sigma_0\left(\dfrac{|x_j|-R_0}{d_{\rm pml}}\right)^m, & R_0<|x_j|<R,
	\end{cases}
	\qquad m\ge 1.
\end{align*}
We use the dimensionless stretching factors $\alpha_j(x_j):=1 + \ii \sigma_j(x_j)$, $j=1,2$, so that the Cartesian uniaxial PML coefficients are
\begin{equation}
	A(x):=
	\begin{pmatrix}
		\alpha_2(x_2)/\alpha_1(x_1) & 0\\
		0 & \alpha_1(x_1)/\alpha_2(x_2)
	\end{pmatrix},
	\qquad
	J(x):=\alpha_1(x_1)\alpha_2(x_2).
	\label{eq:pml-coeffs}
\end{equation}
Consequently, $A=I$ and $J=1$ in $B_{R_0}$, whereas $A$ is complex anisotropic and $J$ is complex-valued in $B_{\rm pml}$.

With these coefficients, the PML-truncated continuous problem is
\begin{equation}
	\begin{aligned}
		-\nabla\cdot(A\nabla u)-\kappa^2Ju &= f
		\quad\text{in }\Omega, \\
		u &= 0\quad\text{on }\Gamma_{\rm out}.
	\end{aligned}\label{eq:pml-strong}
\end{equation}
Here $f$ denotes the right-hand side after the PML transformation and zero extension to $\Omega$. For later use in the discretization and multiscale analysis, we first fix the notation for inner products and norms. For a subdomain $G\subset\Omega$, $(\cdot,\cdot)_G$ denotes the complex $L^2(G)$ inner product with conjugation in the second argument, and let $\norm{\cdot}_G$ be the associated $L^2(G)$ norm. When $G=\Omega$, the subscript is omitted and simply write $(\cdot, \cdot)$ and $\norm{\cdot}$. We also use the $\kappa$-weighted norm $\norm{v}_{1,\kappa,G}^2:=\norm{\nabla v}_G^2+\kappa^2 \norm{v}_G^2$. The trial and test space for the truncated problem is $V:=H_0^1(\Omega)$. The weak formulation of \eqref{eq:pml-strong} is: find $u\in V$ such that
\begin{align}
	a(u,v):=(A\nabla u,\nabla v)-\kappa^2(Ju,v)
	=(f,v)
	\qquad\forall v\in V .
	\label{eq:pml-weak}
\end{align}
In this paper, we employ a polynomial absorption profile extending to the outer truncation boundary, without an outer constant-absorption region. Although Chen, Wu\cite{ChenWu2008}  consider this setting, their truncation error estimate relies on an additional well-posedness assumption for the PML layer problem and the associated stability constant. In contrast, the wavenumber-explicit analyses of Chen, Xiang\cite{ChenXiang2013}  and Wen, Du\cite{WenDu2026} require an outer region with constant absorption. These results therefore do not directly provide all the quantitative estimates needed here. Accordingly, we assume stability of the truncated UPML problem, with constant \(C^{0}_{\mathrm{sta}}\), together with the required PML truncation error bound. Following standard uniaxial PML truncation analyses \cite{Bramble2013,Galkowski2023}, we impose the following standing assumption on the continuous PML approximation.
\begin{assumption}\label{ass:pml-truncation}
	The PML parameters are chosen so that the truncated PML problem \eqref{eq:pml-strong} is stable and its solution agrees with the outgoing solution in the physical region up to an exponentially small truncation error. Namely, there exist constants $C^{0}_{sta}, C_{\rm pml}, c_{\rm pml}>0$, independent of the mesh sizes(maybe depentent on $\kappa$), for which
    \begin{equation*}
        \inf_{u \in V\setminus\{0\}} \sup_{v \in V\setminus\{0\}} \frac{\abs{a(u,v)}}{ \norm{u}_{1,\kappa} \norm{v}_{1,\kappa}} \ge C^{0}_{\rm sta} > 0, 
    \end{equation*}
	\begin{equation*}
		\norm{u^{\rm ext}-u}_{H^1(B_{R_0})} \le C_{\rm pml}\exp(-c_{\rm pml}\kappa d_{\rm pml}) \norm{f^{\rm phys}}_{L^2(B_{R_0})}.
	\end{equation*}
\end{assumption}

For later use in the discretization and multiscale analysis, we first fix the notation for inner products and norms. For a subdomain $G\subset\Omega$, $(\cdot,\cdot)_G$ denotes the complex $L^2(G)$ inner product with conjugation in the second argument, and let $\norm{\cdot}_G$ be the associated $L^2(G)$ norm. When $G=\Omega$, the subscript is omitted and simply write $(\cdot, \cdot)$ and $\norm{\cdot}$. We also use the $\kappa$-weighted norm $\norm{v}_{1,\kappa,G}^2:=\norm{\nabla v}_G^2+\kappa^2 \norm{v}_G^2$. The trial and test space for the truncated problem is $V:=H_0^1(\Omega)$. The weak formulation of \eqref{eq:pml-strong} is: find $u\in V$ such that
\begin{align}
	a(u,v):=(A\nabla u,\nabla v)-\kappa^2(Ju,v)
	=(f,v)
	\qquad\forall v\in V .
	\label{eq:pml-weak}
\end{align}
The rest of the paper treats this PML-truncated variational problem as the continuous model. In the next section, we introduce nested finite element spaces on $\Omega$ and augment the fine-scale PML discretization with continuous interior penalty terms. This yields the CIP-PML reference problem on which the subsequent LOD construction is built.

\section{CIP-PML discretization and preliminary estimates}\label{sec:cip}
\subsection{Meshes, finite element spaces, and patches}\label{sec:mesh}
To discretize the PML problem \eqref{eq:pml-weak}, we follow the standard LOD setting \cite{Gallistl2015,Malqvist2014,Peterseim2017} and introduce a regular and quasi-uniform triangulation $\TT_H$ of $\Omega$ with mesh size $H$. In addition, let $\TT_h$ be a conforming refinement of $\TT_H$, with mesh size $h < H$.  For a simplicial mesh $\TT$ of a two-dimensional domain and a subdomain $G \subset \Omega$, we denote by $\faces_i(\TT, G), i = 0,1,2$, the set of all $i$-dimensional faces of $\TT$ contained in $\overline{G}$. In the special case $G = \Omega$, we simply write $\faces_i(\TT) := \faces_i(\TT, \Omega)$. In particular, $\faces_0(\TT_H)$ and $\faces_1(\TT_H)$ are the sets of coarse vertices and coarse edges of $\TT_H$, respectively. We further denote by $\ifaces_i(\TT)$ the subset of $i$-faces that are not contained in the boundary $\partial \Omega$. For a simplex $K \in \faces_2(\TT)$, let $\mathbb P_1(K)$ denote the space of scalar polynomials of total degree at most one on $K$. We define the conforming fine-scale finite element space by
\begin{equation*}
	V_{h}:=\{v_h \in V : v_h|_k \in \mathbb P_1(k) \; \text{for all } k \in \faces_2(\TT_h) \}.
\end{equation*}
The corresponding fine-grid Galerkin approximation of the PML problem is: find $u_h \in V_h$ such that
\begin{align}
	a_{h}(u_h, v_h) = (f, v_h) \quad \forall v_h \in V_h. \label{eq:fine-weak}
\end{align}
The coarse finite element space is defined analogously on $\TT_H$:
\begin{equation*}
	V_{H}:=\{ v_H \in V : v_H|_K \in \mathbb P_1(K) \; \text{for all } K \in \faces_2(\TT_H)\}. 
\end{equation*}
Since $\TT_h$ is a conforming refinement of $\TT_H$, we have $V_H \subset V_h \subset V$. The space $V_H$ represents the coarse scale in the subsequent multiscale construction. Since homogeneous Dirichlet boundary condition is imposed on $\partial \Omega$, only interior coarse vertices give rise to independent nodal basis functions. Then $V_H$ is spanned by the standard Lagrange basis $(\phi_{z})_{z\in \faces_0(\TT_H)}$.

We shall frequently use coarse patches in the construction and analysis of localized correctors. For a coarse element $T \in \faces_2(\TT_H)$, define the element patch $N^L(T), L \in \mathbb N$ recursively by $N^0(T) := \operatorname{int}(T)$, and, for $L \ge 1$,
\begin{equation*}
	N^L(T) := \operatorname{int} \bigg( \bigcup \bigl\{ K\in \faces_2(\TT_H) : \overline K\cap \overline{N^{L-1}(T)}\neq\emptyset \bigr\} \bigg).
\end{equation*}
Thus $N^L(T)$ is obtained from $T$ by adding $L$ layers of coarse elements. This patch construction satisfies a standard finite-overlap property. More precisely, we define 
\begin{equation}\label{overlapping}
	C_{{\rm ol},L} := \max_{T\in\TT_H} \# \bigl\{ K\in \faces_2(\TT_H) : K\subset \overline{N^L(T)} \bigr\}.
\end{equation}
By the shape regularity and quasi-uniformity of $\TT_H$, there exists a constant $C_{\rm ol}^\# > 0$, independent of $H, h, \kappa$, and $L$, such that $C_{{\rm ol},L} \le C_{\rm ol}^{\#}(1+L)^d$.

\subsection{CIP-stabilized fine-scale PML discretization} 
We now introduce the continuous interior penalty stabilization on the fine-scale PML discretization. The penalty is imposed on the interior edges of the fine mesh, namely on $\ifaces_1(\TT_h)$. For an edge $e = \partial k^+ \cap \partial k^- \in \ifaces_1(\TT_h)$, we denote by $n_{k^+}$ and $n_{k^-}$ on $\partial k^+$ and $\partial k^-$ the unit outward normals, respectively. We further define the jump of the normal derivative by
\begin{equation*}
	\jump{\partial_n v_h} := (\nabla v_h|_{k^+})\cdot n_{k^+} + (\nabla v_h|_{k^-})\cdot n_{k^-}.
\end{equation*}
This definition is independent of the labeling of the two neighboring elements. We write $\duality{\phi}{\psi}_e := \int_e \phi \overline\psi \dd \sigma$ for the complex $L^2(e)$ inner product, with conjugation in the second argument. For a family of complex penalty parameters $\{\gamma_e\}_{e \in \ifaces_1(\TT_h)}$, we define the CIP penalty form by 
\begin{equation*}
	j_h^\gamma(u_h, v_h) := \sum_{e\in \ifaces_1(\TT_h)} \gamma_e h_e \duality{\jump{\partial_n u_h}}{\jump{\partial_n v_h}},
\end{equation*}
where $h_e$ denotes the diameter of the edge $e$. Throughout the paper, the parameters $\gamma_e$ are assumed to be uniformly bounded complex numbers with  nonpositive imaginary parts. We specifically denote the real part of $\gamma_{e}$ by $\gamma_{\alpha}$ and the imaginary part by $\gamma_{i}$.
The CIP-stabilized PML sesquilinear form is then given by
\begin{equation}
	a_h^\gamma(u_h, v_h) = a_h(u_h, v_h) + j_h^\gamma(u_h, v_h) \quad \forall u_h, v_h \in V_h. \label{eq:cip-bilinear}
\end{equation}
Equivalently, 
\begin{equation*}
	a_h^\gamma(u_h, v_h) = (A \nabla u_h, \nabla v_h) - \kappa^2(J u_h, v_h) + \sum_{e\in \ifaces_1} \gamma_e h_e \duality{\jump{\partial_n u_h}}{\jump{\partial_n v_h}}.
\end{equation*}
Then the first order CIP-FEM for \eqref{eq:fine-weak} reads as: find $u_{h}^{\gamma}\in V_{h}$ such that 
\begin{align}
	a_{h}^{\gamma}(u_{h}^{\gamma},v_{h})=(f,v_{h}) \quad \forall v_{h}\in V_{h}. \label{eq:cip-weak}
\end{align}
The penalty term is consistent for sufficiently smooth functions, since the normal derivative has no jump across interior edges. Its role is to modify the discrete phase and stability properties of the finite element operator while preserving the conforming Lagrange finite element space. Readers may refer to \cite{Li2019,Wu2014} for detailed studies on the CIP penalty. 

The analysis below requires elementary bounds for the PML-weighted gradient and mass terms. Since the dependence on the PML strength $\sigma_0$ enters the stability and localization estimates, we record these bounds explicitly in the following lemma.
\begin{lemma}\label{lem:pml-coeff}
	The PML coefficients $A$ and $J$ defined in \eqref{eq:pml-coeffs} satisfy the following estimates. For all $u, v \in V_h$,
	\begin{equation*}
		\abs{(A\nabla u, \nabla v)} \leq  p(\sigma_0) \norm{\nabla u}  \norm{\nabla v},
	\end{equation*}
	and
	\begin{equation*}
		\Real (A\nabla u ,\nabla u) \geq \tilde p(\sigma_0) \norm{\nabla u}^{2},
	\end{equation*}
	where one may take $p(\sigma_0) := \sqrt{1 + \sigma_0^2}$ and $\tilde p(\sigma_0) := (1+\sigma_0^2)^{-1}$. Moreover,
	\begin{equation*}
		\abs{(Ju, v)} \le q(\sigma_0) \norm{u} \norm{v},
	\end{equation*}
	with $q(\sigma_0) := 1+\sigma_0^2$. Finally, the real part of the $J$-weighted mass term satisfies 
	\begin{equation*}
		\norm{u}^2_{B_{R_0}} + (1-\sigma_0^2) \norm{u}^2_{B_{\mathrm pml}} \le \Real(Ju, u) \le \norm{u}^2.
	\end{equation*}
\end{lemma}

\begin{proof}
	Recall that $\alpha_j = 1 + \ii \sigma_j, \; j = 1,2$, and 
	\begin{equation*}
		A = \begin{pmatrix}
			\alpha_2 / \alpha_1 & 0 \\ 0 & \alpha_1 / \alpha_2
		\end{pmatrix}, \quad J = \alpha_1 \alpha_2.
	\end{equation*}
	We first estimate the diffusion coefficient $A$. Since 
	\begin{equation*}
		\abs{\frac{\alpha_2}{\alpha_1}} = \frac{\sqrt{1+\sigma^2_2}}{\sqrt{1+\sigma^2_1}} \le \sqrt{1+\sigma^2_0},
	\end{equation*}
	the Cauchy--Schwarz inequality gives
	\begin{align*}
		\abs{(A\nabla u, \nabla v)} &\le \norm{\nabla u}_{B_{R_0}}  \norm{\nabla v}_{B_{R_0}} + \sqrt{1+\sigma_0^2} \norm{\nabla u}_{B_{\mathrm pml}} \norm{\nabla v}_{B_{\mathrm pml}} \\
		& \le \left( \norm{\nabla u}^2_{B_{R_0}} + \sqrt{1+\sigma_0^2} \norm{\nabla u}^2_{B_{\mathrm pml}} \right)^{1/2} \left( \norm{\nabla v}^2_{B_{R_0}} + \sqrt{1+\sigma_0^2} \norm{\nabla v}^2_{B_{\mathrm pml}} \right)^{1/2} \\
		& \le \sqrt{1+\sigma_0^2} \norm{\nabla u} \norm{\nabla v}.
	\end{align*}
	Hence $p(\sigma_0) = \sqrt{1+\sigma_0^2}$. For the real part, direct expansion yields
	\begin{equation*}
		\Real \frac{\alpha_2}{\alpha_1} = \frac{1+\sigma_1\sigma_2}{1+\sigma_1^2} \ge \frac1{1+\sigma_0^2},
	\end{equation*}
	and similarly $\Real \frac{\alpha_1}{\alpha_2} \ge (1+\sigma_0^2)^{-1}.$
	Therefore, pointwise,
	\begin{equation*}
		\Real (A\xi\cdot \xi) \ge \frac1{1+\sigma_0^2}\abs{\xi}, \quad \xi \in \mathbb{C}^2.
	\end{equation*}
	Taken $\xi = \nabla u$ and integrating over $\Omega$ gives
	\begin{equation*}
		\Real (A \nabla u, \nabla u) \ge \frac1{1+\sigma_0^2}\norm{\nabla u}^2.
	\end{equation*}
	Thus $\tilde p(\sigma_0) = (1+\sigma_0^2)^{-1}$. Next, 
	\begin{equation*}
		\abs{(Ju, v)} \le \abs{J} \norm{u} \norm{v} \le (1+\sigma_0^2) \norm{u} \norm{v},
	\end{equation*}
	which proves the $J$-weighted continuity estimate. Finally, $J = (1+\ii \sigma_1)(1+\ii \sigma_2)$ yields $\Real J = 1 - \sigma_1 \sigma_2$. In $B_{R_0}, \sigma_1 = \sigma_2 = 0$, whereas in $B_{\mathrm pml}, 1-\sigma_0^2 \le \Real J \le 1$. Hence
	\begin{equation*}
		\Real (Ju, u) = \int \Real J \abs{u}^2 \dd x
	\end{equation*}
	satisfies
	\begin{equation*}
		\norm{u}^2_{B_{R_0}} + (1-\sigma_0^2) \norm{u}^2_{B_{\mathrm pml}} \le \Real(Ju, u) \le \norm{u}^2.  \qedhere
	\end{equation*}
\end{proof}

To prepare for the multiscale construction, we next introduce the discrete energy norm associated with the CIP-PML formulation and collect the stability properties of the fine-grid reference problem. These estimates will later be combined with the coarse-fine decomposition to establish well-posedness and decay for the corrector problems.

\begin{definition}
	For any subdomain $G \subset \Omega$, which is a union of fine-grid elements, and $v_h \in V_h$, define 
	\begin{equation*}
		\gsnorm{v_h}{G}^2 := \norm{\nabla v_h}_G^2 + \sum_{e\in \ifaces_1(\TT_h, G)} \abs{\gamma_e} h_e  \norm{\jump{\partial_n v_h}}_e^2,
	\end{equation*}
	and
	\begin{equation*}
		\gnorm{v_h}{G}^2 := \abs{v_h}_{1,\gamma,G}^2 + \kappa^2 \norm{v_h}_G^2.
	\end{equation*}
	When $G = \Omega$, we omit the subscript $G$.
\end{definition}

\begin{lemma}\label{lem:penalty}
	There exists a constant $\alpha_0>0$, such that if $-\alpha_0 \le \Real \gamma_e \le \abs{\gamma}$, then for every $\phi_h \in V_h$, 
	\begin{equation*}
		\Real (A\nabla \phi_h, \nabla \phi_h)_G + \Real \sum_{e\in \ifaces_1(\TT_h, G)} \gamma_e h_e \duality{\jump{\partial_n\phi_h}}{\jump{\partial_n\phi_h}}_e \ge \frac{\tilde p(\sigma_0)}{2} \gsnorm{\phi_h}{G}^2.
	\end{equation*}
\end{lemma}

\begin{proof}
	By Lemma~\ref{lem:pml-coeff}, $\Real (A\nabla\phi_h, \nabla\phi_h)_G \ge \tilde p(\sigma_0)\|\nabla\phi_h\|_{L^2(G)}^2$. Moreover, since $\Real\gamma_e \ge -\alpha_0$, a standard discrete trace/inverse estimate yields
	\begin{equation*}
		\Real\sum_{e\in\ifaces_1(\TT,G)} \gamma_e h_e
		\duality{\jump{\partial_n\phi_h}}{\jump{\partial_n\phi_h}}_e \ge - \alpha_0\sum_{e\in\ifaces_1(\TT,G)} h_e\norm{\jump{\partial_n\phi_h}}_e^2 \ge - \alpha_0 C_{\rm inv} \norm{\nabla \phi_h}_G^2.
	\end{equation*}
	Hence
	\begin{equation*}
		\Real (A\nabla \phi_h, \nabla \phi_h)_G + \Real \sum_{e\in \ifaces_1(\TT,G)} \gamma_e h_e \duality{\jump{\partial_n\phi_h}}{\jump{\partial_n\phi_h}}_e \ge (\tilde p(\sigma_0) - \alpha_0 C_{\rm inv}) \norm{\nabla \phi_h}_G^2.
	\end{equation*}
	On the other hand, by the uniform bound $|\gamma_e| \le \abs{\gamma}$ and the same inverse estimate,
	\[
	\gsnorm{\phi_h}{G} \le (1+C_{\rm inv}|\gamma|)\norm{\nabla\phi_h}_G^2.
	\]
	Therefore, choosing $\alpha_0>0$ sufficiently small such that $\tilde p(\sigma_0)-\alpha_0 C_{\rm inv} \ge \tilde p(\sigma_0)(1+C_{\rm inv}|\gamma|)/2$, the asserted inequality follows, which completes the proof.
\end{proof}

\begin{lemma}
	Assume that the penalty parameters $\gamma_e$ satisfy the condition in \Cref{lem:penalty}. Then the sesquilinear form $a_h^\gamma$ defined in \eqref{eq:cip-bilinear} satisfies the discrete inf-sup condition
	\begin{equation}
		\inf_{u_h \in V_h\setminus\{0\}} \sup_{v_h \in V_h\setminus\{0\}} \frac{\abs{a_h^\gamma(u_h,v_h)}}{\ognorm{u_h} \ognorm{v_h}} \ge C_{\rm sta} > 0, \label{eq:fine-infsup}
	\end{equation}
	where $C_{\rm sta}$ is independent of the particular choice of $u_h$ and $v_h$. Consequently, for every right-hand side $f\in V_h'$, the discrete problem \eqref{eq:cip-weak} admits a unique solution $u_h^\gamma\in V_h$.
\end{lemma}

\begin{proof}
	We define the elliptic projection operators $P^{\pm }_{h}: H^{2}(\Omega)\to V_{h}$ as follows (cf. \cite[Section~4.2]{Li2019}):
	\begin{align}\label{ell proj}
		&\begin{aligned}
			& (AP_{h}^{+}\nabla v, \nabla \phi_h)_G +  \sum_{e\in \ifaces_1(\TT_h, G)} \gamma_e h_e \duality{\jump{\partial_n P^{+}_{h}v}}{\jump{\partial_n\phi_h}}_e \\
			=  &(A\nabla v, \nabla \phi_h)_G +  \sum_{e\in \ifaces_1(\TT_h, G)} \gamma_e h_e \duality{\jump{\partial_nv}}{\jump{\partial_n\phi_h}}_e\quad \forall \phi_{h}\in V_{h},
		\end{aligned}
		\\
		&\begin{aligned}
			 &(A\nabla \phi_{h}, \nabla P_{h}^{-}v)_G +  \sum_{e\in \ifaces_1(\TT_h, G)} \gamma_e h_e \duality{\jump{\partial_n \phi_{h}}}{\jump{\partial_n P^{+}_{h}v}}_e \\
			= &(A\nabla \phi_{h}, \nabla v)_G +  \sum_{e\in \ifaces_1(\TT_h, G)} \gamma_e h_e \duality{\jump{\partial_n \phi_{h}}}{\jump{\partial_n v}}_e \quad \forall \phi_{h}\in V_{h}.
		\end{aligned}
	\end{align}
	Under the assumption on the penalty parameter in \Cref{lem:penalty}, we can obtain the same elliptic projection error. Based on this elliptic projection, we can directly follow the proof of the pre‑asymptotic error estimates in \cite[Theorem~4.4]{Li2019} to derive the stability estimates of the CIP‑FEM for the UPML problem; see \cite{Li2019,Li2025}.
\end{proof}

\subsection{Quasi-interpolation and fine-scale kernel}\label{sec:quasi}
The LOD construction is driven by a coarse quasi-interpolation operator whose kernel represents the unresolved fine scale. Since the PML coefficients $A$ and $J$ are uniformly bounded above and below by \Cref{lem:pml-coeff}, and since the present analysis does not involve high-contrast media, we employ a standard unweighted local $L^2$-averaging quasi-interpolation operator of Cl\'ement type; see, for instance, \cite{Malqvist2014,Malqvist2020,Peterseim2017}. For the weighted adjustment of the quasi-interpolation projection operator underlying the LOD multiscale methods in high-contrast heterogeneous media, as well as for implementation aspects of the corresponding high-contrast LOD, we refer the reader to \cite{peterseim2020computational}. 

For each interior coarse vertex $z \in \ifaces_0(\TT_H)$, let $\phi_z$ denote the associated nodal basis function of $V_H$. We define $\II_H : V_h \to V_H$ by
\begin{equation*}
	\II_H v := \sum_{z \in \ifaces_0(\TT_H)} \alpha_z (v) \phi_z, \qquad \alpha_z(v) := \frac{(v, \phi_z)}{(1, \phi_z)}.
\end{equation*}
Since each basis function $\phi_z$ is supported in the nodal patch $\omega_z := \supp \phi_z$, the coefficient $\alpha_z(v)$ depends only on the values of $v$ in $\omega_z$. The operator $\II_H$ satisfies the standard local stability and approximation estimate: for every $T \in \faces_2(\TT_H)$, 
\begin{equation*}
	H^{-1}\norm{v - \II_H v}_T + \norm{\nabla(v-\II_H v)}_T \le C_{\II_H} \norm{\nabla v}_{N^1(T)}, \quad v \in V_h,
\end{equation*}
where the constant $C_{\II_H}$ depends only on the shape regularity of $\TT_H$; see \cite{clement1975approximation} for the proof of local error estimates of $\mathcal{I}_{H}$. We define the fine-scale space by 
\begin{equation*}
	W_h := \ker(\II_H) = \{w_h \in V_h : \II_H w_h = 0\}.
\end{equation*}
This space collects those fine-grid functions that are invisible on the coarse scale and servers as the detail space in all subsequent corrector problem. Although $\II_H$ is not a projection, it admits a bounded right inverse $\RR_H:V_H \to V_h$ with local support preservation: for every $v_H \in V_H$, 
\begin{equation*} 
	\II_H \RR_H v_H = v_H, \quad \norm{\nabla \RR_H v_H} \le C_{\II_H}' \norm{\nabla v_H}, \quad \supp(\RR_H v_H) \subseteq \supp(v_H),
\end{equation*}
where $C_{\II_H}'$ depends only on the shape regularity of $\TT_H$; see \cite{Engwer2019,Peterseim2017} for the explicit construction.

Consequently, every $v_h \in V_h$ can be decomposed as
\begin{equation*}
	v_h = \RR_H \II_H v_h + \left(v_h - \RR_H \II_H v_h \right), \quad v_h - \RR_H \II_H v_h \in W_h,
\end{equation*}
which provides the basis for separating coarse and fine scales and formulating the ideal and localized corrector problems in the next section.

\section{LOD for the CIP-FEM}\label{sec:lod}
In this section, we construct the ideal and localized LOD schemes for the fine-grid CIP-PML problem \eqref{eq:cip-weak}. The derivation follows the variational multiscale principle and is based on the quasi-interpolation operator introduced in \Cref{sec:quasi}.

\subsection{An Ideal Method} 
Recall that the fine-scale space is $W_h := \ker(\II_H) \subset V_h$. Let $(\mathcal{I}_{H}|_{V_{H}})^{-1}: V_{H}\to V_{H}$ be the bounded inverse for $\mathcal{I}_{H}|_{V_{H}}$, and define $Q_H := (\mathcal{I}_{H}|_{V_{H}})^{-1}\II_H: V_h \to V_H$. Then, for every $v_h \in V_h$, 
\begin{equation*}
	v_h = Q_H v_h + (\II- Q_H)v_h, \quad (\II- Q_H)v_h \in W_h.
\end{equation*}
Indeed, $Q_{H}$ is exactly the $L^{2}$ projection from $V_{h}$ onto $V_{H}$, see \cite{maalqvist2015computation,Peterseim2017}. With this convention, we write 
\begin{align*}
	V_h = V_H + W_h.
\end{align*}

For $v_h \in V_h$, we introduce the conjugate-linear functional $\A(v_h)(\cdot)$ associated with $a_h^\gamma$,
\begin{equation}\label{eq:Agama}
	\A(v_h)(w_h) := a_h^\gamma(v_h, w_h), \quad v_h, w_h \in V_h.
\end{equation}
Similarly, we define the linear functional $\aA(v_h)(\cdot)$ associated with the adjoint form by 
\begin{equation*}
	\aA(v_h)(w_h) := {a_h^\gamma(w_h, v_h)}, \quad v_h, w_h \in V_h.
\end{equation*}

\begin{assumption}\label{asp:resolution}
	(Resolution condition)Given the constant $C_{\II_H}$ from the local stability estimate of $\II_H$, the fine-overlap constant $C_{\rm ol, 1}$ from \Cref{sec:mesh}, and the coefficient $q(\sigma_0)$ from \Cref{lem:pml-coeff}, we assume that the coarse mesh size $H$ satisfies 
	\begin{equation*}
		C^2_{\II_H} C_{\rm ol, 1} q(\sigma_0) (\kappa H)^2 \le \frac{\tilde p(\sigma_0)}{3}.
	\end{equation*}
\end{assumption}

\begin{definition}\label{def: Ideal Corrector}
	For any bounded conjugate-linear functional $\textbf{F} \in W_h^\prime$, the ideal Green's corrector $\green F \in W_h$ is defined by 
	\begin{equation}\label{eq:green}
		a_h^\gamma (\green \textbf{F}, w_h) = \textbf{F}(w_h) \quad \forall w_h \in W_h. 
	\end{equation}
	Analogously, the adjoint Green's corrector $\agreen \textbf{F} \in W_h$ is defined by 
	\begin{equation}\label{eq:adj-green}
		a_h^\gamma (w_h, \agreen \textbf{F}) = \textbf{F}(w_h) \quad \forall w_h \in W_h. 
	\end{equation}
\end{definition}
The next estimate shows that \eqref{eq:green} and \eqref{eq:adj-green} are both well-posed on $W_h$.

\begin{lemma}\label{lem:ker-coer}
	Under \Cref{asp:resolution}, the sesquilinear form $a_h^\gamma$ is coercive on $W_h$. More precisely, for every $w_h \in W_h$,
	\begin{equation*}
		\Real a_h^\gamma(w_h, w_h) \ge \frac{\tilde p(\sigma_0)}{6} \ogsnorm{w_h}^2.
	\end{equation*}
\end{lemma}

\begin{proof}
	Since $w_h \in W_h = \ker(\II_H)$, the local approximation property of $\II_H$ and the finite-overlap argument imply
	\begin{equation*}
		\kappa^2 q(\sigma_0) \norm{w_h}^2 \le C_{\II_H}^2 C_{\rm ol, 1} q(\sigma_0)(\kappa H)^2 \ogsnorm{w_h}^2.
	\end{equation*}
	By \Cref{asp:resolution}, the right-hand side is bounded by $\tilde p(\sigma_0) \ogsnorm{w_h}^2 / 3$. Combining this with \Cref{lem:penalty} yields
	\begin{equation*}
		\Real a_h^\gamma(w_h, w_h) \ge \frac{\tilde p(\sigma_0)}{2} \ogsnorm{w_h}^2 - \frac{\tilde p(\sigma_0)}{3} \ogsnorm{w_h}^2 = \frac{\tilde p(\sigma_0)}{6} \ogsnorm{w_h}^2. \qedhere
	\end{equation*}
\end{proof}

\begin{remark}\label{remark:3UPML}
	(i) The wavenumber-independent coercivity of $a_h^\gamma$ on the fine-scale space $W_h$, as established in \Cref{lem:ker-coer} under \Cref{asp:resolution}, is a key ingredient in the localization analysis of the LOD method. In particular, it ensures that the constants arising in the local energy estimates are independent of $\kappa$, which is crucial for establishing the exponential decay of the ideal element correctors.
	(ii). For UPML problem in $\Omega \subset \mathbb{R}^{3}$ with complex coordinate stretching factors $\alpha_{j}(x_{j})=1+(\zeta+\ii)\sigma_j(x_{j}), j=1,2,3$ (recall \Cref{sec:model}). The matrix $A$ below is the three-dimensional UPML matrix; see \cite{Bramble2013}. If $\zeta \geq 1$, there exist a constant $C_{\sigma_{0},\zeta}$ such that $\Real(\xx^HA\xx)\geq C_{\sigma_{0},\zeta} \xx^H \xx$ for all $\xx \in \mathbb{C}^3$. Therefore, for sufficiently small $\kappa H$, the variational problem in \Cref{def: Ideal Corrector} is coercive with the coercive constant depending only on $\sigma_{0}$ and $\zeta$. If $\zeta=0$, coercivity no longer holds. But from a recent result in \cite{bai2026puresourcetransferdomain}, we have
	\begin{align}
	 \inf_{u_h\in V_h \setminus \{0\}} \sup_{v_h\in V_h \setminus \{0\}}\frac{|(A\nabla u_{h},\nabla v_{h})+(\lambda\kappa^{2}+\Lambda )(u_{h},v_{h})|}{\parallel u_{h}\parallel^{\gamma } \parallel v_{h}\parallel^{\gamma }}\geq \mu(\sigma_{0}),
	\end{align}
	where $\lambda$, $\Lambda$ and $\mu(\sigma_{0})$ depend only on $\sigma_{0}$. In this case, \Cref{def: Ideal Corrector} is well defined by the Lax–Milgram–Babu\v{s}ka theorem under the proper resolution condition like \Cref{asp:resolution}.
\end{remark}

The next result gives the ideal VMS (Variational Multiscale Scheme) decomposition of the fine-grid solution; see also \cite{altmann2021numerical} for the related abstract framework.

\begin{lemma}\label{lem:vms}
	Under \Cref{asp:resolution}, let $u_h^\gamma \in V_h$ be the solution of \eqref{eq:cip-weak}, and set $u_H := Q_H(u_h^\gamma) \in V_H$. Then there holds
	\begin{equation}\label{eq:decomp}
		u_h^\gamma = u_H - (\green \circ \A) u_H + \green \textbf{F}_f, \quad \textbf{F}_f(\cdot) := (f, \cdot).
	\end{equation}
	Moreover, if $f \in H^j(\Omega)$ with $0 \le j \le 2$, we obtain 
	\begin{equation*}
		H^{-1}(C_{\rm ol, 1} C_{\II_H})^{-1/2} \norm{\green \textbf{F}_f} + \ogsnorm{\green \textbf{F}_f} \le \frac{6}{\tilde p(\sigma_0)} H^{j+1} C_{{\rm proj}, j} (C_{\rm ol, 1} C_{\II_H})^{1/2} \abs{f}_{H^j(\Omega)}.
	\end{equation*}
\end{lemma}

\begin{proof}
	Since $u_h^\gamma - u_H \in W_h$, testing \eqref{eq:cip-weak} with an arbitrary $w_h \in W_h$ gives
	\begin{equation*}
		a_h^\gamma(u_h^\gamma - u_H, w_h) = (f, w_h) - a_h^\gamma(u_H, w_h) = a_h^\gamma(\green \textbf{F}_f - (\green \circ \A) u_H, w_h),
	\end{equation*}
	which proves \eqref{eq:decomp}. To estimate the source corrector, we first observe that $\II_H \green \textbf{F}_f = 0$ implies $(v_H, \green \textbf{F}_f) = 0$ for all $v_H \in V_H$. Let $\Pi_H : L^2(\Omega) \to V_H$ denote the $L^2$ projection onto the coarse space. Then we obtain
	\begin{equation*}
		a_h^\gamma(\green \textbf{F}_f, \green \textbf{F}_f) = (f, \green \textbf{F}_f) = (f - \Pi_H f, \green \textbf{F}_f).
	\end{equation*}
	Applying the coercivity estimate from \Cref{lem:ker-coer}, the Cauchy--Schwarz inequality, and the projection estimate
	\begin{equation*}
		\norm{f - \Pi_H f} \le C_{{\rm proj}, j} H^j \abs{f}_{H^j(\Omega)},
	\end{equation*}
	we obtain 
	\begin{equation*}
		\frac{\tilde p(\sigma_0)}{6} \ogsnorm{\green \textbf{F}_f} \le  C_{{\rm proj}, j} H^j \abs{f}_{H^j(\Omega)} \norm{G_\infty^\gamma(\textbf{F}_f)}.
	\end{equation*}
	Finally, the $L^2$-term is controlled by the kernel estimate.
\end{proof}

The next theorem characterizes the corresponding ideal Petrov--Galerkin scheme.

\begin{theorem}\label{VMSinfsup}
	Under the assumptions of \Cref{lem:vms}, the ideal multiscale approximation $\msI{u} := (\II- \green \circ \A) u_H$ in (\ref{eq:decomp}) is characterized as the unique solution of 
	\begin{equation}\label{eq:pg-weak}
		a_h^\gamma \big((\II- \green \circ \A)u_H, (\II- \agreen \circ \aA)v_H \big) = \big(f, (\II- \agreen \circ \aA)v_H \big) \quad \forall v_H \in V_H.
	\end{equation}
	Moreover, the associated ideal trial and test spaces satisfy the inf-sup condition
	\begin{equation*}
		\inf_{u_H \in V_H\setminus\{0\}} \sup_{v_H \in V_H\setminus\{0\}} \frac{\abs{ a_h^\gamma \big( (\II - \green \circ \A)u_H, (\II - \agreen \circ \aA)v_H \big) }}{\ognorm{(\II - \green \circ \A)u_H} \ognorm{(\II - \agreen \circ \aA)v_H}} \ge C_{\rm ms} > 0.
	\end{equation*}
	Here $C_{\rm ms}$ only depends on $C_{\rm sta}$ and $\sigma_{0}$.
\end{theorem}

\begin{proof}
	For every $v_H \in V_H$, the adjoint corrector satisfies 
	\begin{equation*}
		a_h^\gamma \big( w_h, (\II - \agreen \circ \aA) v_H \big) = 0 \quad \forall w_h \in W_h.
	\end{equation*}
	Combining this orthogonality with the decomposition \eqref{eq:decomp} and the fine-grid problem \eqref{eq:cip-weak}, we obtain \eqref{eq:pg-weak}. The inf-sup stability follows from the well-posedness of the corrector problem on $W_h$, together with the isomorphism between $V_H$ and the ideal multiscale trial and test spaces.
\end{proof}

From now on, for convenience, we introduce the shorthand notation
\begin{equation*}
	\IdC := \green \circ \A, \qquad \aIdC := \agreen \circ \aA.
\end{equation*}

\begin{corollary}
	Under the assumption of \Cref{lem:vms}, the ideal Petrov--Galerkin approximation satisfies
	\begin{equation*}
		\ogsnorm{u_h^\gamma - \msI{u}} \leq \frac{6C_{{\rm proj}, j} (C_{\rm ol, 1} C_{\II_H})^{1/2}}{\tilde p(\sigma_0)} H^{j+1} \abs{f}_{H^j(\Omega)}.
	\end{equation*}
	In addition, 
	\begin{equation*}
		\norm{u_h^\gamma - \msI{u}} \le H(C_{{\rm ol}, 1} C_{\II_H})^{1/2} \ogsnorm{u_h^\gamma - \msI{u}}.
	\end{equation*}
\end{corollary}

\begin{proof}
	The first estimate follows immediately from \Cref{lem:vms}. The $L^2$-bound is the kernel estimate applied to $u_h^\gamma - \msI{u} \in W_h$.
\end{proof}
\begin{remark}
	Compared with conventional LOD methods, the LOD method based on the quasi-interpolation operator $\mathcal{I}_{H}$ in this paper further characterizes a higher order of convergence for the error between the idealized multiscale solution and the exact solution when the source term $f$ has regularity higher than $L^{2}$. This property is generally not shared by LOD methods constructed using other quasi-interpolation projections; see \cite[Example~3.10]{altmann2021numerical} for the specific orders of the error of the idealized LOD method with respect to different quasi-interpolation projections. Moreover, this result can be regarded as a refinement and complement to the error analysis and numerical experiments in \cite{Peterseim2017}.
\end{remark}
\subsection{Element correctors and exponential decay}\label{sec:decay}
For $v_H \in V_H$, the ideal corrector $\IdC v_H$ is decomposed into element contributions,
\begin{equation*}
	\IdC v_H = \sum_{T\in \faces_2(\TT_H)} \IdtC v_H.
\end{equation*}
To this end, let $\At v_H \in W_h^\prime$ be the local functionals such that
\begin{equation*}
	(\A v_H)|_{W_h} = \sum_{T\in \faces_2(\TT_H)} \At v_H,
\end{equation*}
where $\At v_H$ depends only on the values of $v_H$ in $T \cup N^1(T)$. For instance, one may split the volume terms elementwise and distribute each fine-edge penalty contribution equally among the adjacent coarse elements. Accordingly we define
\begin{equation*}
	(\At v_H)(w_h) := (A\nabla v_H, \nabla w_h)_T - \kappa^2 (J v_H, w_h)_T +  \sum_{e \in \ifaces_1(\TT_h, \partial T)} \gamma_e h_e \duality{\partial_{n_T} v_H}{\jump{\partial_{n_T} w_h}}.
\end{equation*}
The corresponding ideal element corrector is then defined by 
\begin{equation*}
	\IdtC v_H := \green \circ \At v_H \in W_h.
\end{equation*}
The dual correctors are defined analogously.

\begin{theorem}\label{thm:decay}
	If the resolution condition of \Cref{asp:resolution} is satisfied, then there exist constants $C_{\rm dec} > 0, \delta > 0$, and $\tau = p(\sigma_0) / \tilde p(\sigma_0)$ such that for every $v_H \in V_H$, $T \in \TT_H$, and $L \in \mathbb N$, the element corrector $\phi := \IdtC v_H$ satisfies 
	\begin{equation*}
		\gsnorm{\phi}{\Omega \setminus N^L(T)}^2 \le C_{\rm dec} \tau^2 e^{-\delta L/\tau} \gnorm{v_H}{N^1(T)}^2.
	\end{equation*}
\end{theorem}

The proof is postponed until after a cutoff construction and an auxiliary kernel-correction lemma.
\begin{definition}\label{def:cutoff}
	For $T \in \faces_2(\TT_H)$ and integers $m < M$, let $\eta_T^{m,M}$ be the continuous, piecewise affine function with respect to $\TT_H$, not necessarily satisfying the homogeneous Dirichlet boundary condition, determined by the nodal values
	\begin{equation*}
		\begin{cases*}
			\eta_T^{m,M}(x) = 0, \quad x \in \faces_0(\TT_H, N^m(T)), \\
			\eta_T^{m,M}(x) = 1, \quad x \in \faces_0(\TT_H, \Omega \setminus N^M(T)), \\
			\eta_T^{m,M}(x) = j(M-m)^{-1}, \quad x \in \faces_0(\TT_H, \partial N^{m+j}(T)),\quad j = 0,1,2,\dots,M-m.
		\end{cases*}
	\end{equation*}
\end{definition}

Then $\eta_T^{m,M}$ satisfies 
\begin{equation*}
	\eta_T^{m,M} = 0 \text{ in } N^m(T), \quad \eta_T^{m,M} = 1 \text{ in } \Omega \setminus N^M(T),
\end{equation*}
and
\begin{equation*}
	\norm{H \nabla\eta_T^{m,M}}_{L^\infty(\Omega)} \le \frac{\beta}{M-m},
\end{equation*}
where $\beta > 0$ depends only on the shape regularity of $\TT_H$. The same cutoff function structure is also found in \cite{Malqvist2014}.

\begin{lemma}\label{lem:bh}
	Let $\II_h$ denote the nodal interpolation operator onto $V_h$. For any $w_h \in W_h$, let $\eta := \eta_T^{L-4, L-3}$ be the cutoff function from the \Cref{def:cutoff}, and define $z_h := \II_h(\eta w_h)$. Then there exists $b_h \in V_h$, given by $b_h = \RR_HI_H z_h$, such that $z_h - b_h \in W_h$ and $\supp(b_h) \subset N^{L-2}(T) \setminus N^{L-5}(T)$. Moreover,
	\begin{align*}
		& h^{-1} \norm{\eta w_h - z_h} + \ogsnorm{\eta w_h - z_h} + H^{-1} \norm{b_h} + \ogsnorm{b_h} \\
		\le \; &(C_{\II_h} C_{\rm inv} + C_b) \big( C_{\II_H} \sqrt{C_{\rm ol, 1}} \norm{H \nabla \eta}_{L^\infty(\Omega)} + 1 \big) \gsnorm{w_h}{N^{L-2}(T) \setminus N^{L-5}(T)},
	\end{align*}
	where the constants depend only on the shape regularity of the meshes, the interpolation constants, and the uniform bound on the penalty parameters.
\end{lemma}

\begin{proof}
	Set $R_L := N^{L-2}(T) \setminus N^{L-5}(T)$. Since $\eta = 0$ in $N^{L-4}(T)$ and $\eta = 1$ outside $N^{L-3}(T)$, the function $\eta w_h$ coincide either with $0$ or with $w_h$ away from the transition layer. Hence the interpolation error $\eta w_h - \II_h(\eta w_h)$ is supported in $R_L$. By the local interpolation estimate, together with standard inverse and trace estimates on the fine mesh,
	\begin{equation*}
		\norm{\eta w_h - z_h} + \ogsnorm{\eta w_h - z_h} \le C_{\II_h} C_{\rm inv} \norm{\nabla(\eta w_h)}_{R_L}.
	\end{equation*}
	Using the product rule,
	\begin{equation*}
		\norm{\nabla(\eta w_h)}_{R_L} \le \norm{\nabla w_h}_{R_L} + \norm{w_h \nabla \eta}_{R_L}.
	\end{equation*}
	The first term is bounded by $\gsnorm{w_h}{R_L}$. For the second term, since $w_h \in W_h = \ker(\II_H)$, the local approximation estimate for $\II_H$, summed over the coarse elements in the transition layer, gives
	\begin{equation*}
		H^{-1} \norm{w_h}_{N^{L-3}(T) \setminus N^{L-4}(T)}  \le C_{\II_H} \sqrt{C_{\rm ol, 1}} \gsnorm{w_h}{R_L}.
	\end{equation*}	
	Therefore,
	\begin{align*}
		\norm{w_h \nabla \eta}_{R_L} \le \norm{H \nabla \eta}_{L^\infty(\Omega)} H^{-1} \norm{w_h}_{R_L}\le C_{\II_H} \sqrt{C_{\rm ol, 1}}\norm{H \nabla \eta}_{L^\infty(\Omega)} \gsnorm{w_h}{R_L}.
	\end{align*}
	It remains to verify the role of $b_h$. By definition, $b_h = \RR_H \II_H z_h$. Since $\RR_H$ is a right inverse of $\II_H$, $\II_H b_h = \II_H \RR_H \II_h z_h = \II_H z_h$. Consequently, $\II_H (z_h - b_h) = 0$, and hence $z_h - b_h \in W_h$. The support inclusion follows from the locality of $\II_H$ and support-preserving property of $\RR_H$. Finally, the stability of $\RR_H$, and the preceding bound for $\nabla(\eta w_h)$ imply
	\begin{equation*}
		H^{-1} \norm{b_h} + \ogsnorm{b_h} \le C_b \big( C_{\II_H} \sqrt{C_{\rm ol, 1}} \norm{H \nabla \eta}_{L^\infty(\Omega)} + 1 \big) \gsnorm{w_h}{R_L}. \qedhere
	\end{equation*}
\end{proof}

\begin{proof}[Proof of \Cref{thm:decay}]
	Throughout the proof, $X \lesssim Y$ means that $X \le CY$, where $C$ is independent of $h, H, \kappa, L, T$, and $v_H$. We write $E_\ell := \Omega \setminus N^\ell(T)$ and $R_L := N^L(T) \setminus N^{L-7}(T)$. It is enough to prove the estimate for $L \ge 7$; the remaining finitely many values of $L$ are absorbed into the generic constant. Let $\phi := \IdtC v_H$. Since the element functional $\At v_H$ is supported in $N^1(T)$, the definition of the ideal element corrector gives the orthogonality relation
	\begin{equation*}
		a_h^\gamma(\phi, w_h) = 0 \quad \forall w_h \in W_h \text{ with } \supp(w_h) \subset \Omega \setminus N^1(T).
	\end{equation*}
	
	We first note the global stability estimate. Testing the corrector equation with $\phi$, using coercivity of $a_h^\gamma$ on $W_h$, the local continuity of $\At$, and the resolution condition in \Cref{asp:resolution}, we obtain
	\begin{equation*}
		\tilde p(\sigma_0) \ogsnorm{\phi}^2 \lesssim \abs{\At v_H (\phi)} \lesssim p(\sigma_0) \gnorm{v_H}{N^1(T)} \ogsnorm{\phi}.
	\end{equation*}
	Hence $\ogsnorm{\phi} \lesssim \tau \gnorm{v_H}{N^1(T)}$.
	
	We next derive the Caccioppoli-type estimate. Recall that
	\begin{equation*}
		d_h^\gamma(u,v) := (A\nabla u, \nabla v) + j_h^\gamma(u, v), \quad a_h^\gamma(u, v) = d_h^\gamma(u, v) - \kappa^2(Ju,v).
	\end{equation*}
	Let $\eta := \eta_T^{L-4, L-3}$ be the cutoff from \Cref{def:cutoff} and set $z_h := \II_h(\eta \phi)$. By \Cref{lem:bh}, there exists $b_h := \RR_H \II_H z_h$ such that $z_h - b_h \in W_h$ and $\supp(b_h) \subset R_L$. Moreover,
	\begin{equation*}
		\ogsnorm{\eta \phi - z_h} + \ogsnorm{b_h} \lesssim \gsnorm{\phi}{R_L}.
	\end{equation*}
	Since the local functional $\At v_H$ vanishes for test functions supported away from $N^1(T)$, and since $z_h - b_h$ is supported outside $N^1(T)$, we have $0 = a_h^\gamma(\phi, z_h - b_h) = d_h^\gamma(\phi, z_h - b_h) - \kappa^2(J\phi, z_h - b_h)$. Using $\eta = 1$ on $E_{L-3}$ and \Cref{lem:penalty}, 
	we obtain 
	\begin{equation*}
		\tilde p(\sigma_0) \gsnorm{\phi}{E_L}^2 \lesssim \Real d_{h, E_L}^\gamma (\phi, \phi) = \Real d_h^\gamma (\phi, \eta\phi) - \RR_\eta
	\end{equation*}
	where $\RR_\eta := \Real d_h^\gamma (\phi, \eta\phi)- \Real d_{h, E_L}^\gamma (\phi, \phi)$. Since $\eta = 0$ in $N^{L-4}(T)$, $\eta \phi = 0$ and $\nabla(\eta \phi) = 0$, we have
	\begin{equation*}
		\big( A \nabla \phi, \nabla(\eta\phi) \big) - \big( A \nabla \phi, \nabla(\eta \phi) \big)_{E_L} = \big( A \nabla\phi, \nabla(\eta\phi) \big)_{N^L(T) \setminus N^{L-4}(T)},
	\end{equation*}
	and the jump term in $\RR_\eta$ is controlled by
	\begin{equation*}
		\sum_{e \in \faces_1(\TT_h, R_L)} \gamma_e h_e \duality{\jump{\partial_n \phi}}{\eta\jump{\partial_n\phi}}_e + \sum_{e \in \faces_1(\TT_h, R_L)} \gamma_e h_e \duality{\jump{\partial_n \phi}}{\jump{\phi\partial_n\eta}}_e \lesssim \gsnorm{\phi}{R_L}^2.
	\end{equation*}
	Therefore 
	\begin{equation}\label{eq:R-eta}
		\begin{aligned}
			\tilde p(\sigma_0) \gsnorm{\phi}{E_L}^2 &\lesssim \abs{d_h^\gamma(\phi, z_h - b_h)} +\abs{d_h^\gamma(\phi, \eta\phi - z_h)} + \abs{d_h^\gamma(\phi, b_h)} + \abs{\RR_\eta} \\
			&\lesssim \abs{\kappa^2(J\phi, z_h - b_h)} + p(\sigma_0) \gsnorm{\phi}{R_L}^2.
		\end{aligned}
	\end{equation}
	It remains to estimate the mass term. Since $\phi, z_h - b_h \in W_h$, the kernel estimate for $\II_H$, together with \Cref{lem:bh}, yields
	\begin{equation*}
		\norm{\phi}_{E_{L-5}} \lesssim H \gsnorm{\phi}{E_{L-6}} \text{ and } \norm{z_h - b_h}_{E_{L-5}} \lesssim H\big( \gsnorm{\phi}{E_{L-6}} + \gsnorm{\phi}{R_L} \big).
	\end{equation*}
	Therefore, there holds
	\begin{equation}\label{eq:mass-term}
		\abs{\kappa^2 (J\phi, z_h - b_h)} \lesssim q(\sigma_0) (\kappa H)^2 \big( \gsnorm{\phi}{E_{L-6}}^2 + \gsnorm{\phi}{R_L}^2 \big).
	\end{equation}
	By \Cref{asp:resolution} and the inclusion $E_{L-6} \subset E_L \cup R_L$, the first term on the right-hand side of \eqref{eq:mass-term} can be absorbed into the the left-hand side of \eqref{eq:R-eta}. Consequently,
	\begin{equation*}
		\tilde p(\sigma_0) \gsnorm{\phi}{E_L}^2 \le C_* p(\sigma_0)\gsnorm{\phi}{R_L}^2.
	\end{equation*}
	Here $C_*$ depends only on the shape regularity of the meshes, the interpolation constants, the stability constant of the right inverse $\RR_H$, the finite-overlap constants of a fixed number of coarse layers, the uniform bounds on the penalty parameters, and the constants in the standard trace and inverse estimates. In particular, $C_*$ is independent of $h,H,\kappa,L,T$, and $v_H$.
	
	Since $\gsnorm{\phi}{R_L}^2 \le \gsnorm{\phi}{E_{L-7}}^2 - \gsnorm{\phi}{E_L}^2$, we arrive at 
	\begin{equation*}
		\gsnorm{\phi}{E_L}^2 \le \theta_\sigma \gsnorm{\phi}{E_{L-7}}^2, \quad \theta_\sigma := \frac{C_* p(\sigma_0)}{\tilde p(\sigma_0) + C_* p(\sigma_0)} = \frac{C_* \tau}{1+C_*\tau} < 1.
	\end{equation*}
	Iteration gives $\gsnorm{\phi}{E_L}^2 \le \theta_\sigma^{\lfloor L/7 \rfloor} \ogsnorm{\phi}^2$. Moreover, $\theta_\sigma = 1 - 1/(1+C_* \tau) \le e^{-c/\tau}$ with a constant $c>0$ independent of $h, H,\kappa,L,T$, and $v_H$. Hence
	\begin{equation*}
		\gsnorm{\phi}{E_L}^2 \lesssim e^{-\delta L / \tau} \ogsnorm{\phi}^2 \le C_{\rm dec} \tau^2  e^{-\delta L / \tau} \gnorm{v_H}{N^1(T)}^2,
	\end{equation*}
	The constants $\delta$ and $C_{\rm dec}$ have the same structural dependence as $C_*$, and are independent of $h,H,\kappa,L,T$, and $v_H$.
\end{proof}

\subsection{Localized correctors}\label{sec:lod-loc-cor}
We now introduce the localized version of the ideal correctors. Throughout this subsection, the integer $L \in \mathbb{N}$ denotes the oversampling parameter. For a coarse element $T \in \faces_2(\TT_H)$, let $N^L(T)$ be the $L$-layer patch defined in \Cref{sec:mesh}. We use the following localized spaces:
\begin{equation*}
	V_h(N^L(T)) := \{ v_h \in V_h : v_h = 0 \in \Omega \setminus N^L(T) \},
\end{equation*}
and 
\begin{equation*}
	\pker := W_h \cap V_h(N^L(T)).
\end{equation*}
When the patch touches the outer boundary, the homogeneous Dirichlet condition is inherited from $V_h$. For $v_H \in V_H$, the localized element corrector $\PtC v_H \in \pker$ is defined by
\begin{equation*}
	a_h^\gamma(\PtC v_h, w_h) = (\At v_H)(w_h) \quad \forall w_h \in \pker.\label{local correctors}
\end{equation*}
Here $\At v_H$ is the local functional introduced in \Cref{sec:decay}. Since $\pker \in W_h$, the coercivity estimate of \Cref{lem:ker-coer} also hold on $\pker$. Hence the localized corrector problem is well posed. The global localized corrector is then obtained by summing all element correctors:
\begin{equation*}
	\PC v_H := \sum_{T \in \faces_2(\TT_H)} \PtC v_H.
\end{equation*}
The adjoint localized element corrector $\aPtC v_H$ is defined analogously using the adjoint variational form, and we set 
\begin{equation*}
	\aPC v_H := \sum_{T \in \faces_2(\TT_H)} \aPtC v_H.
\end{equation*}
Consequently, the localized multiscale trial and test spaces are
\begin{equation*}
	\pmsI{V} := (\II - \PC) V_H, \quad \apmsI{V} := (\II - \aPC) v_H.
\end{equation*}
The local Petrov--Galerkin LOD approximation is then given by: find $\pmsI{u} \in \pmsI{V}$ such that 
\begin{equation}\label{sol: gamma LOD}
	a_h^\gamma(\pmsI{u}, \apmsI{v}) = (f, \apmsI{v}) \quad \forall \apmsI{v} \in \apmsI{V}.
\end{equation}
We next prove that the localized correctors converge exponentially to their ideal counterparts as the oversampling parameter increases. The proof uses a cutoff argument similar to that in \Cref{thm:decay}.
\begin{lemma}\label{lem:bh-local}
	Let $L \ge 4, T \in \faces_2(\TT_H)$, and let $\xi_T^L \in V_H$ be a coarse cutoff function satisfying $\xi_T^L = 1$ in $N^{L-3}(T)$, $\xi_T^L = 0$ in $\Omega \setminus N^{L-2}(T)$, and $\norm{H \nabla \xi_T^L}_{L^\infty(\Omega)} \le C_\xi$, where $C_\xi$ depends only on the shaper regularity of $\TT_H$. Then, for every $w_h \in W_h$, with $z_h := \II_h(\xi_T^L w_h)$, there exists $b_h \in V_h$ such that $z_h - b_h \in \pker$, $\supp(b_h) \subset N^{L-1}(T) \setminus N^{L-4}(T)$, and
	\begin{equation*}
		h^{-1}\norm{\xi_T^L w_h - z_h} + \ogsnorm{\xi_T^L w_h - z_h} + H^{-1} \norm{b_h} + \ogsnorm{b_h} \le C_{\rm cut} \gsnorm{w_h}{N^{L-1}(T) \setminus N^{L-4}(T)}.
	\end{equation*}
	The constant $C_{\rm cut}$ depends only on the shape regularity of the meshes, the stability constant of $\II_H$ and $\II_h$, the inverse and trace constants, and the uniform bound on the penalty parameters.
\end{lemma}

\begin{proof}
	The proof is identical to that of \Cref{lem:bh}, with the cutoff function $\eta_T^{L-4, L-3}$ replaced by $\xi_T^L$. The only difference is the location and orientation of the transition layer. Since $\xi_T^L$ is constant outside $N^{L-2}(T) \setminus N^{L-3}(T)$, the interpolation error $\xi_T^L w_h - \II_h(\xi_T^L w_h)$ is supported in a uniformly enlarged annulus contained in $N^{L-1}(T) \setminus N^{L-4}(T)$. Using the local interpolation estimate for $\II_h$, the product rule, the inverse and trace inequalities on the fine mesh, and the kernel estimate for $w_h \in W_h = \ker \II_H$, we obtain
	\begin{equation*}
		h^{-1}\norm{\xi_T^L w_h - \II_h(\xi_T^L w_h)} + \ogsnorm{\xi_T^L w_h - \II_h(\xi_T^L w_h)} \le C_{\rm cut} \gsnorm{w_h}{N^{L-1}(T) \setminus N^{L-4}(T)}.
	\end{equation*}
	The constant depends only on the shape regularity of the meshes, the interpolation constants, the inverse and trace constants, and the uniform bound on the penalty parameters.
\end{proof}

We now estimate the truncation error of a single element corrector. 
\begin{theorem}\label{thm:el-corr-appr}
	Assume that the resolution condition in \Cref{asp:resolution} holds. Then, for every $T \in \faces_2(\TT_H)$, $v_H \in V_H$,
	\begin{equation}\label{eq:el-corr-appr-1}
		\ognorm{\IdtC v_H - \PtC v_H}^2 \le C_{\rm ele}e^{-\delta L/\tau} \gnorm{v_H}{N^1(T)}^{2}.
	\end{equation}
	The same estimate holds for the adjoint correctors: 
	\begin{equation}\label{eq:el-corr-appr-2}
		\ognorm{\aIdtC v_H - \aPtC v_H}^2 \le C_{\rm ele}e^{-\delta L/\tau} \gnorm{v_H}{N^1(T)}^{2}.
	\end{equation}
	The constant $C_{\rm ele}:=C_{a}C_{\rm cut} C_{\rm dec}\tau^{4}e^{6\delta/\tau}$, $C_{a}$ is an independent constant, and other quantities see \Cref{thm:decay} for details.
\end{theorem}

\begin{proof}
	Let $\phi_T := \IdtC v_H$ and $\phi_{T,L} := \PtC v_h$. Since the localized corrector is the Galerkin projection of the ideal corrector onto $\pker$, we have 
	\begin{equation*}
		a_h^\gamma(\phi_T - \phi_{T,L}, w_h) = 0 \quad \forall w_h \in \pker.
	\end{equation*}
	The coercivity on $W_h$, the continuity of $a_h^\gamma$, and C\'ea's lemma yield \begin{equation}\label{eq:el-corr-appr-3}
		\ognorm{\phi_T - \phi_{T,L}} \lesssim \frac{p(\sigma_0)}{\tilde p(\sigma_0)} \inf_{w_L \in \pker} \ognorm{\phi_T - w_L}.
	\end{equation}
	Let $\xi_T^L$ be the cutoff function from \Cref{lem:bh-local} and choose $w_L := \II_h(\xi_T^L \phi_T) - b_h$, where $b_h := \RR_H \II_H \II_h (\xi_T^L \phi_T)$. By \Cref{lem:bh-local}, $w_L \in \pker$. Using the decomposition 
	\begin{equation*}
		\phi_T - w_L = (\II - \xi_T^L)\phi_T + (\xi_T^L \phi_T - \II_h(\xi_T^L \phi_T)) + b_h
	\end{equation*}
	and applying \Cref{lem:bh-local} and \Cref{thm:decay}, we obtain
	\begin{equation}\label{eq:el-corr-appr-4}
		\ognorm{\phi_T - w_L}^2 \le C_{\rm cut} \gnorm{\phi_T}{\Omega \setminus N^{L-6}(T)}^2 \le C_{\rm cut} C_{\rm dec}\tau^{2}e^{6\delta/\tau}e^{-\delta L/\tau}\gnorm{v_H}{N^1(T)}^2.
	\end{equation}
	Combining \eqref{eq:el-corr-appr-3} and \eqref{eq:el-corr-appr-4}, we now prove \eqref{eq:el-corr-appr-1}. The proof of \eqref{eq:el-corr-appr-2} is identical, with the adjoint form replacing $a_h^\gamma(\cdot,\cdot)$.
\end{proof}

The next results gives the corresponding global corrector estimate.

\begin{theorem}\label{thm:corr-appr}
	Under the assumption of \Cref{thm:el-corr-appr}, there exists a constant $C_{\rm glob} > 0$, independent of $h, H, \kappa$, and $L$, such that
	\begin{equation*}
		\ognorm{(\IdC - \PC) v_H}^2 \le C_{\rm glob} e^{-\delta L/\tau} \ognorm{v_H}^2 \quad \forall v_H \in V_H.
	\end{equation*}
	The same estimate holds for the adjoint corrector:
	\begin{equation*}
		\ognorm{(\aIdC - \aPC) v_H}^2 \le C_{\rm glob} e^{-\delta L/\tau} \ognorm{v_H}^2 \quad \forall v_H \in V_H.
	\end{equation*}
	Here $C_{\rm glob}:=\widetilde C_{\rm glob}\tau^{2}(1+L)^{d}C_{ol,1}C_{\rm ele}$ and $\widetilde{C}_{glob}$ are independent constants; see $C_{ol,1}$ in \ref{overlapping} and $C_{\rm ele}$ in \Cref{thm:el-corr-appr}.
\end{theorem}

\begin{proof}
	We prove the estimate for the primal corrector. The adjoint estimate follows in the same way. For $T \in \faces_2(\TT_H)$, set $e_T := (\IdtC - \PtC) v_H$ and $Z := (\IdC - \PC) v_H = \sum_{T \in \faces_2(\TT_H)} e_T$. Since $Z \in W_h$, the kernel estimate and the coercivity of $a_h^\gamma$ on $W_h$ imply
	\begin{equation*}
		\ognorm{Z}^2 \lesssim \frac{1}{p(\widetilde{\sigma_{0}})}\Real a_h^\gamma(Z, Z) = \frac{1}{p(\widetilde{\sigma_{0}})}\Real \sum_{T \in \faces_2(\TT_H)} a_h^\gamma (e_T, Z).
	\end{equation*}
	For each $T$, choose a coarse cutoff function $\zeta_T$ such that $\zeta_T = 0$ in $N^{L+2}(T)$, $\zeta_T = 1$ in $\Omega \setminus N^{L+3}(T)$, and $\norm{H \nabla \zeta_T}_{L^\infty(\Omega)} \lesssim 1$. Define $z_T := \II_h(\zeta_T Z)$, $b_T := \RR_H \II_H z_T$, and $\chi_T := z_T - b_T$. Then $\chi_T \in W_h$. Moreover, since $z_T = 0$ in $N^{L+2}(T)$, while $Z - z_t$ is supported in $N^{L+4}(T)$, the locality of $\II_H$, the identity $\II_H Z = 0$, and the support-preserving property of $\RR_H$ imply $\supp(b_T) \subset N^{L+4}(T) \setminus N^{L+1}(T)$. Consequently, $\chi_T = 0$ in $N^{L+1}(T)$, and the standard cutoff estimate used in \Cref{lem:bh-local} gives
	\begin{equation*}
		\ognorm{Z - \chi_T} \lesssim \gnorm{Z}{N^{L+4}(T)}.
	\end{equation*}
	Using the locality of the corrector equations, we have \begin{equation*}
		a_h^\gamma(\IdtC v_H, \chi_T) = (\At v_H)(\chi_T) = 0 \text{ and } a_h^\gamma(\PtC v_H, \chi_T) = 0.
	\end{equation*} 
	We further obtain
	\begin{equation*}
		a_h^\gamma(e_T, \chi_T) = 0 \text{ and } a_h^\gamma(e_T, Z) = a_h^\gamma(e_T, Z - \chi_T).
	\end{equation*}
	Therefore,
	\begin{equation}\label{eq:corr-appr-1}
		\begin{aligned}
			\ognorm{Z}^2 &\lesssim \tau\sum_{T \in \faces_2(\TT_H)} \ognorm{e_T} \ognorm{Z-\chi_T} \lesssim \tau\sum_{T \in \faces_2(\TT_H)} \ognorm{e_T} \gnorm{Z}{N^{L+4}(T)} \\
			&\le \tau\bigg( \sum_{T \in \faces_2(\TT_H)}\ognorm{e_T}^2 \bigg)^{1/2} \bigg( \sum_{T \in \faces_2(\TT_H)} \gnorm{Z}{N^{L+4}(T)}^2 \bigg)^{1/2}.
		\end{aligned}
	\end{equation}
	By the finite-overlap property of the patches, 
	\begin{equation}\label{eq:corr-appr-2}
		\sum_{T \in \faces_2(\TT_H)} \gnorm{Z}{N^{L+4}(T)}^{2} \lesssim (1+L)^d \ognorm{Z}^2.
	\end{equation}
	It remains to estimate the right-hand side. By \Cref{thm:el-corr-appr},
	\begin{equation*}
		\ognorm{e_T}^2\le C_{\rm ele}e^{-\delta L/\tau} \gnorm{v_H}{N^1(T)}^{2}.
	\end{equation*}
	Summary over $T$ and using the finite-overlap property of the one-layer patches $N^1(T)$, we get
	\begin{equation}\label{eq:corr-appr-3}
		\sum_{T \in \faces_2(\TT_H)} \ognorm{e_T}^2 \le C_{ol,1}C_{\rm ele}e^{-\delta L/\tau} \ognorm{v_H}^2.
	\end{equation}
	The proof is complete by combining \eqref{eq:corr-appr-1}--\eqref{eq:corr-appr-3}.
\end{proof}
\begin{remark}
	Using \Cref{remark:3UPML} (ii) and similar arguments, the conclusions of \Cref{thm:decay},\Cref{thm:el-corr-appr} and \Cref{thm:corr-appr} can be carried over in parallel to the three-dimensional UPML setting.
\end{remark}

\section{A stability-transfer framework for localized multiscale spaces}\label{sec:framework}
The localized corrector analyzed in the previous section are constructed with the same CIP-PML sesquilinear form as the fine-scale reference problem. This gives a direct localization theory for the standard element-based localized spaces. In the next section, however, we shall consider computationally cheaper localized spaces, e.g., spaces obtained from nodal patches or from local corrector problems with the CIP penalty omitted. Such spaces are no longer covered directly by the element-corrector argument of \Cref{sec:lod-loc-cor}.

The purpose of this section is to formulate a stability-transfer principle on the Hilbert space $V$. We restrict our attention here to a broad class of multiscale methods—variational multiscale methods. Throughout this section, all multiscale spaces are subspaces of a Hilbert space $V$. As noted in \cite{altmann2021numerical}, different variational multiscale methods are often induced by different choices of the fine-scale space $W$:
\begin{align}
    V=V_{H}\oplus W.
\end{align}
We denote by $b: V \times V \to \mathbb C$ the sesquilinear form used in the Petrov--Galerkin paring. To emphasize the nature of the stability of the stability-transfer argument, we use the notation $\tnorm{\cdot}$ for the norm, in which continuity, stability, and localization errors are measured. In the present CIP-LOD application, $V=V_{h}$, $b = a_h^\gamma$, and the norm $\tnorm{\cdot}$ is precisely the discrete energy norm $\tnorm{v_h} := \ognorm{v_h}$ for all $v_h \in V_h$.

We assume that $b$ is continuous with respect to $\tnorm{\cdot}$, namely
\begin{equation*}
	\abs{b(u, v)} \le C_b \tnorm{u} \tnorm{v} \quad \forall u, v \in V.
\end{equation*}
\subsection{Idealized variational multiscale method}\label{VMSSS}
For convenience, we first define the multiscale trial and test spaces in the VMS framework in \cite[Section~6.2]{altmann2021numerical}:
\begin{equation*}
	\msII{V}:=\{\msII{u}\in V: b(\msII{u}, w) = 0, \quad\forall w \in W\},
\end{equation*}
and the test spaces defined as:
\begin{equation*}
	\amsII{V}:=\{\amsII{v}\in V: b(w, \amsII{v}) = 0, \quad \forall w \in W\}.
\end{equation*}
Under the stability of the space decomposition and the well-posedness of the original problem, the well-posedness of the idealized variational multiscale method is usually readily obtained; for a detailed proof, see the discussion of \cref{VMSinfsup} or \cite{peterseim2020computational}. Here we state the well-posedness of the idealized variational multiscale method directly as an assumption:
\begin{assumption}\label{asp:abs-infsup}
	There exists a constant $\msII{\alpha} > 0$ such that
	\begin{equation*}
		\inf_{\msII{u} \in \msII{V}\setminus\{0\}} \sup_{\amsII{v} \in \amsII{V}\setminus\{0\}} \frac{\abs{b(\msII{u}, \amsII{v})}}{\tnorm{\msII{u}} \tnorm{\amsII{v}}} \ge \msII{\alpha}.
	\end{equation*}
\end{assumption}
As noted above, different choices of \(W\) induce different multiscale variational formulations. Reference \cite[Section~6]{altmann2021numerical} lists several numerical schemes, such as the local Green's function approach \cite{hughes1998}, residual-free bubbles \cite{brezzi1996}, and classical MsFEM \cite{hou1997} which have the same $W:=\oplus_{T\in \mathcal{T}_{H} }H_{0}^{1}(T)$. The multiscale basis functions constructed by such methods already have local compact support and therefore require no additional localization. Note that localization here is meant relative to the entire space, rather than the oversampling idea introduced in MsFEM to handle oscillations across coarse-element boundaries.
\begin{definition}\label{QOIS}
    Macroscopic quantities of interest (QoIs). Let $N\in\mathbb{N}$. We introduce a family of continuous coarse-scale separating functionals $q_i \in V^\prime, i = 1,\dots, N$. Moreover, we define a so-called information operator $\PP v := \big( q_1(v), \dots, q_N(v) \big) \in \mathbb C^N$.
\end{definition}
Beyond the above choices for $W$, previous studies took a more abstract view and constructed the space \(W\) using quantities of interest in Definition \ref{QOIS}. \cite{altmann2021numerical} gives several examples of specific choices for QoIs. Building on this, more commonly used VMS methods typically choose \(W:=\ker \PP\) via approaches such as LOD \cite{altmann2021numerical, Malqvist2020, Peterseim2017}, operator-adapted wavelets \cite{owhadi1945}, and GMsFEM combined with spectral multiscale techniques \cite{cheung2018}.
\begin{remark}
    In LOD, QoIs are typically encoded implicitly through a quasi-interpolation operator \(\pi_H: V \to V_{\mathrm{aux}}\), see \cite{Peterseim2017,peterseim2020computational,verfürth2017}; hence \(\{q_i\}_{i\in N}\) can be represented as \(\{\pi_H^* \phi_i^*\}_{i\in N}\). Here, \(\pi_{H}^{*}\) denotes the adjoint of \(\pi_{H}\), and \(\{\phi_i^*\}_{i\in N}\) is a family of dual basis functions for the coarse auxiliary space \(V_{\mathrm{aux}}\). The same holds for CEM-MsFEM; see \cite[Section~3]{cheung2018}.
\end{remark}
The multiscale basis functions constructed via QoIs are often globally supported and therefore must be localized in practical computations. Depending on the variational problem and the choice of QoIs, a wide variety of localization strategies may arise. For variational multiscale methods requiring localization, the well-posedness analysis and error estimates often demand extensive modifications when the localization strategy changes, and may even rely on specialized argument techniques. Moreover, in existing well-posedness and error estimates based on classical localization strategies, the assumptions on localization are frequently not the weakest possible \cite{Peterseim2017,peterseim2020computational}, and some error estimates are not optimal \cite{verfürth2017}. Our analysis indicates that these issues are largely due to overestimation arising from the structure of the arguments. This section aims to establish a transfer operator used solely for analysis, which links the idealized space to the localized space, thereby yielding a unified abstract framework and streamlining the well-posedness and error estimate arguments.

Related abstract stability results have been developed in the context of variational multiscale methods; see, for instance, \cite{altmann2021numerical}. Stability and error estimates are derived within a framework that is tied to a particular VMS localization construction. Our approaches is formulated differently. It does not require an explicit representation of the localization procedure. Instead, by employing the transfer operator, we demonstrate that, under certain assumptions on the transfer operator, the well-posedness and error estimates for VMS problems with various localization strategies—such as LOD, SLOD \cite{Freese2024}, and MS-GFEM \cite{MaAlberScheichl2023} can be established within a unified framework. Importantly, the assumptions imposed here are no stronger than those typically required in the VMS framework of \cite{altmann2021numerical}. More significantly, when the transfer operator preserves QoIs, the analysis of well-posedness and error estimates becomes more complete than existing arguments. Specifically, the oversampling conditions can be relaxed in theory, the proof structure becomes unified, and for special problems the conclusions can be formally optimal. The formulation is carried out on an arbitrary Hilbert space $V$. It is therefore not restricted to scalar $H^1$-conforming problems. This abstract setting also provides a basis for extending the stability-transfer principle to variational multiscale methods in $\bf H(\rm curl)$ and $\bf H(\rm div)$.
\subsection{Transfer operators}
In this section, we first introduce an abstract stability-transfer operator that links the local multiscale function spaces to the idealized multiscale space. Section \ref{1} employs this transfer operator to derive stability and error estimates for localized variational multiscale methods, thereby unifying the stability and error analyses for different localization strategies. Section \ref{VMS2} shows that if the transfer operator possesses certain separating functional properties, the abstract transfer-operator argument can improve existing proofs and yield results under weaker conditions and with clearer conclusions.

\begin{assumption}\label{asp:abs-decay}
There exist bijective transfer operators \(\mathcal{Q}_{L}^{ms}:V_{\infty}^{ms}\to V_{L}^{ms}\) and \(\mathcal{Q}_{L}^{ms,\star}:V_{\infty}^{ms,\star}\to V_{L}^{ms,\star}\) whose associated transfer errors satisfy the following estimates: there exists a quantity $\varepsilon_L > 0$ such that
	\begin{equation*}
		\tnorm{(\II - \Q) \msII{\phi}} \le \varepsilon_L \tnorm{\msII{\phi}} \quad \forall \msII{\phi} \in \msII{V},
	\end{equation*}
	and
	\begin{equation*}
		\tnorm{(\II - \aQ) \amsII{\phi}} \le \varepsilon_L \tnorm{\amsII{\phi}} \quad \forall \amsII{\phi} \in \amsII{V}.
	\end{equation*}
\end{assumption}

For standard localized LOD spaces, one typically has $\varepsilon_L =  C_{\rm glob} e^{-\delta L/\tau}$ (recall \Cref{thm:corr-appr}). The following lemma records a direct consequence of Assumption \ref{asp:abs-decay}.

\begin{lemma}\label{lem:abs-decay}
	If Assumption \ref{asp:abs-decay} holds and $\varepsilon_L < 1$, then $\Q$ and $\aQ$ are isomorphisms. Moreover, for all $\msII{\phi} \in \msII{V}$,
	\begin{equation*}
		(1 - \varepsilon_L) \tnorm{\msII{\phi}} \le \tnorm{\Q \msII{\phi}} \le (1 + \varepsilon_L)\tnorm{\msII{\phi}},
	\end{equation*}
	and the same estimate holds for $\aQ$. Equivalently, if $\pmsII{\phi} = \Q \msII{\phi}$, then
	\begin{equation*}
		\frac1{1+\varepsilon_L} \tnorm{\pmsII{\phi}} \le \tnorm{\msII{\phi}} \le \frac1{1-\varepsilon_L}\pmsII{\phi}.
	\end{equation*}
\end{lemma}

\begin{proof}
	The upper bound follows from the triangle inequality:
	\begin{equation*}
		\tnorm{\Q \msII{\phi}} \le \tnorm{\msII{\phi}} + \tnorm{(\II - \Q) \msII{\phi}} \le (1 + \varepsilon_L) \tnorm{\msII{\phi}}.
	\end{equation*}
	The lower bound follows similarly:
	\begin{equation*}
		\tnorm{\msII{\phi}} \ge \tnorm{\msII{\phi}} - \tnorm{(\II - \Q) \msII{\phi}} \ge (1 - \varepsilon_L)\tnorm{\msII{\phi}}.
	\end{equation*}
	The proof for the adjoint transfer operator is identical.
\end{proof}
\subsubsection{Well-posedness and error estimates via the abstract transfer operator}\label{1}
\begin{lemma}\label{VMSINFSUPLOCAL}
    Under Assumption \ref{asp:abs-decay}, and we suppose $2C_{b}\varepsilon_{L}(2+\varepsilon_{L})<\alpha_{\infty}^{ms}$. Then the localize Petrov-Galerkin pair satisfies
    \begin{equation*}
        \inf_{\pmsII{u} \in \pmsII{V}\setminus\{0\}} \sup_{\apmsII{v} \in \apmsII{V}\setminus\{0\}} \frac{\abs{b(\pmsII{u}, \apmsII{v})}}{\tnorm{\pmsII{u}} \tnorm{\apmsII{v}}}=\alpha_{L}^{ms} \ge \frac{\alpha_{\infty}^{ms}-2C_{b}\varepsilon_{L}(2+\varepsilon_{L})}{2(1+\varepsilon_{l})^{2}}.
    \end{equation*}
\end{lemma}
Let $F \in V^\prime$ be a bounded right-hand side. The ideal multiscale solution $\msII{\BB}F \in \msII{V}$ is defined by 
\begin{equation}\label{vmsidealsolu}
	b(\msII{\BB} F, \amsII{v}) = F(\amsII{v}) \quad \forall \amsII{v} \in \amsII{V}.
\end{equation}
The localized multiscale solution $\pmsII{\BB}F \in \pmsII{V}$ is defined by 
\begin{equation}\label{vmslocalsolu}
	b(\pmsII{\BB}F, \apmsII{v}) = F(\apmsII{v}) \quad \forall \apmsII{v} \in \apmsII{V}.
\end{equation}
\begin{lemma}\label{VMSERROR}
    Under assumptions of \cref{VMSINFSUPLOCAL} and $\varepsilon_{L}<1$, we have 
    \begin{align}
        \tnorm{\msII{\BB} F-\pmsII{\BB}F}\le(1+\frac{1}{\alpha_{L}^{ms}})\varepsilon_{L}\tnorm{\msII{\BB}F}+\frac{\varepsilon_{L}}{\alpha_{L}^{ms}(1-\varepsilon_{L})}(\parallel F\parallel_{V_{h}^{'}}+C_{b}\tnorm{\msII{\BB}F}).
    \end{align}
\end{lemma}
\subsubsection{Transfer operator preserving QoIs}\label{VMS2}
Recalling the fine-scale space \(W_{\PP}\) introduced in \cref{VMSSS}, we consider the VMS induced by \(W_{\PP}\). Furthermore, if the transfer operator is assumed to preserve the coarse-scale degrees of freedom $\{q_{i}\}_{i=1,\dots,N}$, the theoretical proof becomes more complete.
\begin{assumption}\label{asp:abs-unisol}
	The map $\PP$ is unisolvent on the ideal and localized multiscale spaces. More precisely, $\PP |_{\msII{V}} : \msII{V} \to \mathbb C^N$ and $\PP |_{\pmsII{V}} : \pmsII{V} \to \mathbb C^N$ are both isomorphisms. The same assumption is imposed on the adjoint spaces:  $\PP |_{\amsII{V}} : \amsII{V} \to \mathbb C^N$ and $\PP |_{\apmsII{V}} : \apmsII{V} \to \mathbb C^N$.
\end{assumption}
\begin{definition}\label{asp:transfer-oper}
	The trial transfer operator $\Q : \msII{V} \to \pmsII{V}$ is defined as follows: for every $\msII{\phi} \in \msII{V}, \Q \msII{\phi} \in \pmsII{V}$ is the unique function satisfying
	\begin{equation*}
		\PP(\Q \msII{\phi}) = \PP(\msII{\phi}).
	\end{equation*}
	Similarly, the adjoint transfer operator $\aQ : \amsII{V} \to \apmsII{V}$ is defined by 
	\begin{equation*}
		\PP(\aQ \amsII{\phi}) = \PP(\amsII{\phi}) \quad \forall \amsII{\phi} \in \amsII{V}.
	\end{equation*}
	The transfer errors belong to the fine-scale kernel:
	\begin{equation*}
		(\II - \Q)\msII{\phi} \in W_\PP, \; \text{ and } \; (\II - \aQ)\amsII{\phi} \in W_\PP.
	\end{equation*}
\end{definition}
In the VMS framework in \cite[Section~3.4]{altmann2021numerical}, \Cref{asp:abs-unisol} can be verified through saddle point problem.
We now impose the abstract approximation property of the localized spaces.

We now prove that sufficiently accurate localized spaces inherit this inf-sup stability.

\begin{theorem}\label{thm:abs-loc-infsup}
	Let Assumption~\ref{asp:abs-unisol}, \ref{asp:abs-decay} and \ref{asp:abs-infsup} hold. Suppose that $C_b \varepsilon_L^2 < \msII{\alpha}$. Then the localized Petrov-Galerkin pair satisfies
	\begin{equation*}
		\inf_{\pmsII{u} \in \pmsII{V}\setminus\{0\}} \sup_{\apmsII{v} \in \apmsII{V}\setminus\{0\}} \frac{\abs{b(\pmsII{u}, \apmsII{v})}}{\tnorm{\pmsII{u}} \tnorm{\apmsII{v}}} \ge \pmsII{\alpha},
	\end{equation*}
	where $\pmsII{\alpha} := (\msII{\alpha} - C_b \varepsilon_L^2)/(1+\varepsilon_L)^2 > 0$. In particular, if $C_b \varepsilon_L^2 \le \msII{\alpha}/2$, then $\pmsII{\alpha} \ge \msII{\alpha}/2(1+\varepsilon_L)^2$.
\end{theorem}

\begin{proof}
	For every $0 \neq \pmsII{u} \in \pmsII{V}$, by \Cref{lem:abs-decay}, there exists a unique $\msII{u} \in \msII{V}$ such that $\pmsII{u} = \Q \msII{u}$. By the ideal inf-sup condition, for any $\delta > 0$, there exists $\amsII{v} \in \amsII{V}$ such that 
	\begin{equation*}
		\abs{b(\msII{u}, \amsII{v})} \ge (\msII{\alpha} - \delta) \tnorm{\msII{u}} \tnorm{\amsII{v}}.
	\end{equation*}
	Set $\apmsII{v} := \aQ \amsII{v}$. Let $e_u := (\II - \Q)\msII{u}$ and $e_v := (\II - \aQ) \amsII{v}$. Then $e_u, e_v \in W_\PP$, $\pmsII{u} = \msII{u} - e_u$ and $\apmsII{v} = \amsII{v} - e_v$. Using the orthogonality relations in the defination of $\msII{V}$ and $\amsII{V}$, we obtain
	\begin{equation*}
		b(\pmsII{u}, \apmsII{v}) = b(\msII{u}, \amsII{v}) + b(e_u, e_v).
	\end{equation*}
	Therefore, by the continuity of $b$ and \Cref{asp:abs-decay},
	\begin{equation*}
		\abs{	b(\pmsII{u}, \apmsI{v})} \ge (\msII{\alpha} - \delta) \tnorm{\msII{u}} \tnorm{\amsII{v}} - C_b \varepsilon_L^2 \tnorm{\msII{u}} \tnorm{\amsII{v}}.
	\end{equation*}
	Combining with \Cref{lem:abs-decay}, and then letting $\delta \to 0$, yields
	\begin{equation*}
		\sup_{\apmsII{v} \in \apmsII{V}\setminus\{0\}} \frac{\abs{b(\pmsII{u}, \apmsII{v})}}{\tnorm{\apmsII{v}}} \ge \frac{\msII{\alpha} - C_b \varepsilon_L^2}{(1+\varepsilon_L)^2} \tnorm{\pmsII{u}}. \qedhere
	\end{equation*}
\end{proof}
The localized problem is well posed whenever the condition of \Cref{thm:abs-loc-infsup} is satisfied. For the error estimate, define the fine-scale residual of $F$ relative to the separating functionals by 
\begin{equation*}
	\eta_\PP(F) := \inf_{\bold a} \norm{F - \sum_{i=1}^{N} a_i q_i}_{W_{\PP}^\prime},
\end{equation*}
where for a given $G \in W_\PP^\prime$,
\begin{equation*}
	\norm{G}_{W_{\PP}^\prime} := \sup_{w \in W_{\PP} \setminus \{0\}} \frac{\abs{G(w)}}{\tnorm{w}}.
\end{equation*}
This quantity measures the part of the right-hand side that is invisible to the coarse separating functionals.

\begin{theorem}\label{thm:err}
	Let the assumptions in \Cref{thm:abs-loc-infsup} hold. If $\varepsilon_L < 1$, then
	\begin{equation*}
		\tnorm{\msII{\BB}F - \pmsII{\BB}F} \le \bigg( \varepsilon_L + \frac{C_b \Lambda_L}{\pmsII{\alpha}} \varepsilon_L^2 \bigg) \tnorm{\msII{\BB}F} + \frac{\Lambda_L}{\pmsII{\alpha}} \varepsilon_L \eta_\PP(F),
	\end{equation*}
	where $\Lambda_L := (1-\varepsilon_L)^{-1}$. Consequently, if $\Lambda_{L}\varepsilon_L\leq \alpha_L^{ms}$, then
	\begin{equation*}
		\tnorm{\msII{\BB}F - \pmsII{\BB}F} = \mathcal O(\varepsilon_L)\tnorm{\msII{\BB}F} +  \eta_\PP(F).
	\end{equation*}
\end{theorem}

\begin{proof}
	We split the error as 
	\begin{equation*}
		\msII{\BB}F - \pmsII{\BB}F = (\msII{\BB}F - \Q \msII{\BB}F) + (\Q\msII{\BB}F - \pmsII{\BB}F).
	\end{equation*}
	The first term is bounded by \Cref{asp:abs-decay}. It remains to estimate $\pmsII{E} := \Q\msII{\BB}F - \pmsII{\BB}F \in \pmsII{V}$. By the localized inf-sup condition,
	\begin{equation*}
		\pmsII{\alpha} \tnorm{\pmsII{E}} \le \sup_{\apmsII{v} \in \apmsII{V}\setminus\{0\}} \frac{\abs{b(\pmsII{E}, \apmsII{v})}}{\tnorm{\apmsII{v}}}.
	\end{equation*}
	Let $\amsII{v} := (\aQ)^{-1} \apmsII{v} \in \amsII{V}$. Set $e_u := (\II - \Q) \msII{\BB}F$ and $e_v := (\II - \aQ)\amsII{v}$. Then $e_u, e_v \in W_{\PP}$, $\Q \msII{\BB}F = \msII{\BB}F - e_u$, and $\apmsII{v} = \amsII{v} - e_v$. Using the definition of $\msII{\BB}$ and $\pmsII{\BB}$, and the orthogonality relations, we obtain
	\begin{align*}
		b(\pmsII{E}, \apmsII{v}) & = b(\Q\msII{\BB}F, \apmsII{v}) - F(\apmsII{v}) \\ 
		&= b(\msII{\BB}F-e_u, \amsII{v} - e_v) - F(\amsII{v} - e_v) \\
		&= F(\amsII{v}) + b(e_u, e_v) - F(\amsII{v}) + F(e_v) \\
		&= F(e_v) + b(e_u, e_v).
	\end{align*}
	Since $e_v \in W_{\PP}$, for any $\bold a = (a_1, \dots, a_N) \in \mathbb C^N$,
	\begin{equation*}
		\sum_{i=1}^N a_i q_i(e_v) = 0.
	\end{equation*}
	Subsequently,
	\begin{equation*}
		\abs{F(e_v)} = \abs{\bigg( F - \sum_{i=1}^N a_i q_i \bigg)(e_v)} \le \eta_\PP(F) \tnorm{e_v}.
	\end{equation*}
	By \Cref{asp:abs-decay} and \Cref{lem:abs-decay},
	\begin{equation*}
		\tnorm{e_v} \le \varepsilon_L\tnorm{\amsII{v}} \le \Lambda_L \varepsilon_L \tnorm{\apmsII{v}}.
	\end{equation*}
	Similarly, $\tnorm{e_u} \le \varepsilon_L \tnorm{\msII{\BB}F}$. Combining with the continuity of $b$, we get
	\begin{equation*}
		\abs{b(\pmsII{E}, \apmsII{v})} \le \Lambda_L \varepsilon_L \eta_\PP(F) \tnorm{\amsII{v}} + C_b \Lambda_L \varepsilon_L^2 \tnorm{\msII{\BB}F} \tnorm{\apmsII{v}}.
	\end{equation*}
	Taking the supremum over $\apmsII{v}$ gives
	\begin{equation*}
		\tnorm{\pmsII{E}} \le \frac{\Lambda_L}{\pmsII{\alpha}} \varepsilon_L \eta_\PP(F) + \frac{C_b \Lambda_L}{\pmsII{\alpha}} \varepsilon_L^2 \tnorm{\msII{\BB}F}. \qedhere
	\end{equation*}
\end{proof}

\begin{remark}
	The abstract residual $\eta_\PP(F)$ vanishes if the right-hand side is represented exactly by the separating functionals on $W_\PP$. In general, $\eta_\PP(F)$ measures the fine-scale component of the forcing that cannot be captured by the coarse degrees of freedom. This is consistent with the source corrector term in the ideal LOD decomposition, see \Cref{lem:vms}.
\end{remark}
\begin{remark}
	These abstract results can apply directly to more complex localization strategies, provided that Assumptions \ref{asp:abs-unisol}, \ref{asp:abs-decay}, and \ref{asp:abs-infsup} can be verified for the corresponding localized spaces and transfer operators, e.g., Super-Localized Orthogonal Decomposition \cite{Freese2024}. Once the required approximation and preservation properties are established, Theorems \ref{thm:abs-loc-infsup} and \ref{thm:err} apply directly. This avoids repeating a method-specific stability analysis. It also yields inf-sup stability and localization-error estimates in a form consistent with those of the classical LOD method.
\end{remark}

For the CIP-LOD method of the present paper, the abstract spaces may be chosen as $\msII{V} = (\II - \IdC)V_H$ and $\amsII{V} = (\II - \aIdC)V_H$, and similarly for the localized spaces (see \Cref{local correctors}). The functionals $q_i$ are taken as the coarse degrees of freedom induced by $\II_H$. Since the correctors belong to the fine-scale kernel, the transfer error lies in $W = \ker\II_H$. Most importantly, in the CIP‑LOD method, the abstractly defined transfer operator in \Cref{asp:transfer-oper} takes the explicit form $(\mathcal{I}-\mathcal C_{L}^{\gamma})Q_{H}$; see \Cref{local correctors} for details. In the specific CIP-LOD method considered in \Cref{sec:lod}, the remaining mathematical quantities appearing in the abstract framework of this section have explicit expressions; see \Cref{sec:lod} for details.

Under the above settings, we can naturally obtain the error approximation property of the CIP-LOD solution
\begin{corollary}
	Under the resolution assumption in \Cref{asp:resolution} and overlapping condition 
	\begin{equation*}
		C_{a^{\gamma}} C_{\rm glob} e^{-\delta L/\tau}\le C_{\rm ms}/2,
	\end{equation*}
	where $C_{a^{\gamma}}$ is the norm of $a_{h}^{\gamma}$ and see \Cref{thm:corr-appr} for $C_{\rm glob}$. The error between the localized LOD approximation $u_{H,L,ms}^{\gamma}\in V_{H,L,ms}^{\gamma}$ in (\ref{sol: gamma LOD}) and the fine-scale background solution $u_{h}$ in (\ref{eq:cip-weak}) can be expressed as follows:
	\begin{equation*}
		\norm{u_h - \pmsI{u}}_\gamma \leq C_{\rm err, 1}(\sigma_0)  C_{\rm glob}^{1/2} e^{-\delta L/2\tau} \norm{f} + C_{\rm err, 2}(\sigma_0)H^{j+1}\abs{f}_{H^j(\Omega)}, \quad 0 \le j \le 2.
	\end{equation*}
	Here $C_{\rm err, 1}$ and $C_{\rm err, 2}$ depend only on $\sigma_{0}$. 
\end{corollary}
\section{Improved CIP-LOD method}\label{sec:alg-opt}

The localized correctors in \Cref{sec:lod} are constructed with the same CIP-PML from $a_h^\gamma$ as the fine-scale reference problem. This is natural for analysis, but it requires the assembly of CIP edge contributions in every local corrector problem. In this section we describe a cheaper variant: the local correctors are computed without the CIP penalty, whereas the final Petrov--Galerkin problem is still assembled with the full stabilized form $a_h^\gamma$.

This localized trial and test spaces are $\pmsIII{V}{0} := (\II - C^{0}_{L})V_H$ and $\apmsIII{V}{0} := (\II - C^{0,\star}_{L})V_H$. And the corresponding Petrov--Galerkin method is: find $\pmsIII{u}{0} \in \pmsIII{V}{0}$ such that
\begin{equation}\label{sol:without CIP}
	a_h^\gamma(\pmsIII{u}{0}, \apmsIII{v}{0}) = (f, \apmsIII{v}{0}) \quad \forall \apmsIII{v}{0} \in \apmsIII{V}{0}.
\end{equation}
According to the abstract framework of \Cref{sec:framework}, we observe that the space $\pmsIII{V}{0}$ can be regarded as a certain localization strategy for the space $\pmsIII{V}{\gamma}$, and the transfer operator $\mathcal{Q}_{L}^{ms}$ in \Cref{asp:transfer-oper} can be explicitly expressed as:
\begin{align*}
	\mathcal{Q}_{L}^{ms}:=(\mathcal{I}-\mathcal{C}_{L}^{0})Q_{H}.
\end{align*}
In order to correctly utilize the abstract framework, we only need to show that the approximation error of the transfer operator satisfies \Cref{asp:abs-decay}.


\begin{theorem}\label{thm:err without CIP}
	Assume \Cref{asp:resolution} and the overlapping condition
	\begin{align*}
		2 C_b C^{2}(\kappa H,\sigma_{0},L)(\kappa^{2} h^{2}+e^{-\delta L/\tau})\le \frac{C_{\rm ms}}{2}
	\end{align*} 
	hold. The error between the solution $u_{H,L,ms}^{\gamma}$ in (\ref{sol: gamma LOD}) and $u^{0}_{H,L,ms}$ in (\ref{sol:without CIP}) can be estimated by
	\begin{align}
		\norm{\pmsI{u} - \pmsIII{u}{0}}_\gamma \le C(\kappa H, \sigma_0, L) (\kappa h + e^{-\delta L/2\tau}) \norm{f},
	\end{align}
	where $C(\kappa H, \sigma_0, L)$ is independent of $h$ and $v_H$.
\end{theorem}

Before proving \Cref{thm:err without CIP}, we establish several auxiliary results for comparing the global CIP and non-CIP multiscale basis functions. These results will be used to estimate the perturbation caused by omitting the CIP stabilization in the local corrector problems. 

To this end, we first introduce a saddle-point characterization of the global multiscale basis functions. Let $\mathcal{A}: V \to V^\prime, \duality{\mathcal{A} u}{v} := a(u, v)$, denote the operator associated with the sesquilinear form of the truncated PML problem; see \ref{eq:pml-weak}. Let $\varphi_z^{\rm ms} := (\II - \mathcal{C}_\infty)\varphi_z$ denote the corresponding global multiscale basis function. By construction,it preserves the coarse-scale component $Q_H \varphi_z^{\rm ms} = \varphi_z$, and is orthogonal to the fine-scale space, namely $a(\varphi_z^{\rm ms}, w) = 0$ for all $w \in W$. Hence, $\mathcal{A} \varphi_z^{\rm ms}$ lies in the range of $Q_H^T$, which is the adjoint operator of $Q_H$. Therefore, there exists a Lagrange multiplier $\lambda_{\varphi_z} \in V_H$ such that $\mathcal{A} \varphi_z^{\rm ms} + Q_H^T \lambda_{\varphi_z} = 0$. Together with the coarse-scale constraint, this yields the saddle-point system
\begin{equation*}
    \setlength{\arraycolsep}{2pt} 
    \begin{pmatrix}
        \mathcal{A} & Q_H^T \\
        Q_H & 0
    \end{pmatrix}
    \begin{pmatrix}
        \varphi_z^{\rm ms} \\
        \lambda_{\varphi_z}
    \end{pmatrix}
    =
    \begin{pmatrix}
        0 \\
        \varphi_z
    \end{pmatrix}.
\end{equation*}
The same characterization applies to the CIP and non-CIP global multiscale basis function. More precisely, replacing $\mathcal{A}$ by $\mathcal{A}^\gamma$ introduced in \ref{eq:Agama}, or by $\mathcal{A}^0: V_h \to V_h^\prime, \duality{\mathcal{A}^0 u}{v} := a_h^0(u,v),$ gives the saddle-point formulations for $\varphi_z^{\rm ms,\gamma} := (\II - \IdC) \varphi_z$ and $\varphi_z^{\rm ms, 0} := (\II - \mathbb{C}_\infty^0 \varphi_z)$, respectively.  

We next define the coarse-scale operator $\mathcal{S}: V_H \to V_H$ by the relation $S \lambda_{\varphi_z} := \varphi_z$. The well-definedness and bijectivity of $\mathcal{S}$ will be established in \Cref{lem:schur} below. If $\mathcal{A}$ is invertible, the first block equation gives $\varphi_z^{\rm ms} = -\mathcal{A}^{-1} Q_H^T \lambda_H$. Substituting this expression into $Q_H \varphi_z^{\rm ms} = \varphi_z$ yields $- Q_H \mathcal{A}^{-1} Q_H^T \lambda_{\varphi_z} = \varphi_z$. Hence, whenever $\mathcal{A}$ is invertible, $\mathcal{S}$ admits the Schur-complement representation
\begin{equation*}
	\mathcal{S} := -Q_H \mathcal{A}^{-1}Q_H^T.
\end{equation*}

\begin{lemma}\label{lem:schur}
	The coarse-scale operator $\mathcal{S}$ is an isomorphism. Moreover, for every $\varphi_z \in V_H$,
	\begin{equation*}
		\norm{\mathcal{S}^{-1}\varphi_{z}} \lesssim (H^{-2} + \kappa^2) \norm{\varphi_{z}}
	\end{equation*}
	The corresponding ideal multiscale function $\varphi_z^{\rm ms}$ satisfies
	\begin{equation}\label{eq:schur}
		a(\varphi_z^{\rm ms}, v) = -(\mathcal{S}^{-1}\varphi_z, v) \quad \forall v \in V.
	\end{equation}
\end{lemma}
\begin{proof}
	We first show that $\mathcal{S}$ is injective. Suppose that $\varphi_z = \mathcal{S} \lambda_{\varphi_z} = 0$. Then $Q_H \varphi_z^{\rm ms} = 0$, so that $\varphi_z^{\rm ms} \in W$. Moreover, testing the first block equation with an arbitrary $w \in W$ yields
	\begin{equation*}
		0 = a(\varphi_z^{\rm ms}, w) + \duality{Q_H^T \lambda_{\varphi_z}}{w} = a(\varphi_z^{\rm ms}, w) + (\lambda_{\varphi_z}, Q_H w) = a(\varphi_z^{\rm ms}, w).
	\end{equation*}
	 The coercivity of $a$ on $W$ therefore implies $\varphi_z^{\rm ms} = 0$. The first block equation then reduces to $Q_H^T \lambda_{\varphi_z} = 0$. Since $Q_H$ is surjective, $Q_H^T$ is injective, and hence $\lambda_{\varphi_z} = 0$. Hence, $\mathcal{S}$ is injective. Since $V_H$ is finite dimensional, $\mathcal{S}$ is an isomorphism. The saddle-point characterization then gives
	 \begin{equation*}
	 	a(\varphi_z^{\rm ms}, v) = -(\lambda_{\varphi_z}, v) = -(\mathcal{S}^{-1}\varphi_z, v) \quad \forall v \in V.
	 \end{equation*}
 	Furthermore, the stability of the ideal corrector and the inverse estimate on $V_H$ imply
 	\begin{equation*}
 		\norm{\lambda_{\varphi_z}}^2 = \abs{a(\varphi_z^{\rm ms}, \lambda_{\varphi_z})} \lesssim \norm{\varphi_z^{\rm ms}}_{1, \kappa} \norm{\lambda_{\varphi_z}}_{1, \kappa}  \lesssim (H^{-2} + \kappa^2) \norm{\varphi_z} \norm{\lambda_{\varphi_z}}.
 	\end{equation*}
 	Consequently, $\norm{\mathcal{S}^{-1} \varphi_z} \lesssim (H^{-2} + \kappa^2) \norm{\varphi_z}$.
\end{proof}

\begin{lemma}\label{lem:H2}
	Let \Cref{asp:resolution} hold. Then the ideal global multiscale basis function $\varphi_z^{\rm ms}$ satisfies
	\begin{equation*}
		 \kappa \norm{\varphi_z^{\rm ms}} + \norm{\nabla \varphi_z^{\rm ms}} \lesssim 1.
	\end{equation*} 
	Moreover,
	\begin{equation*}
		\abs{\varphi_z^{\rm ms}}_{H^2} \lesssim (1 + \frac1{\kappa H})\kappa.
	\end{equation*}
	The hidden constants are independent of $h$, $H$, and $\kappa$.
\end{lemma}

\begin{proof}
	The inf-sup stability of the ideal corrector problem, the continuity of $a$, the local support of $\varphi_z$, and \Cref{asp:resolution} imply 
	\begin{equation*}
		\kappa \norm{\varphi_z^{\rm ms}} + \norm{\nabla \varphi_z^{\rm ms}} \lesssim \norm{\varphi_z^{\rm ms}}_{1, \kappa} \lesssim \norm{\varphi_z}_{1,\kappa} \lesssim \sqrt{1 + (\kappa H)^2} \lesssim 1.
	\end{equation*}
	Recalling the definition of $a$ and (\ref{eq:schur}), we obtain
	\begin{equation*}
		(A\nabla \varphi_z^{\rm ms}, \nabla v) = (\kappa^2 J \varphi_z^{\rm ms} - \mathcal{S}^{-1} \varphi_z, v) \quad \forall v \in V.
	\end{equation*}
	Standard $H^2$-regularity for the principle PML operator therefore yields 
	\begin{equation*}
		\abs{\varphi_z^{\rm ms}}_{H^2} \lesssim \kappa^2 \norm{\varphi_z^{\rm ms}} + \norm{\mathcal{S}^{-1} \varphi_z}.
	\end{equation*}
	The first estimate gives $\kappa^2 \norm{\varphi_z^{\rm ms}} \lesssim \kappa$. Furthermore, \Cref{lem:schur}, $\norm{\varphi_z}\lesssim H$, and $\kappa H \lesssim 1$ imply
	\begin{equation*}
		\abs{\varphi_z^{\rm ms}}_{H^2} \lesssim \kappa + (H^{-2} + \kappa^2)\norm{\varphi_z} \lesssim (1 + \frac1{\kappa H})\kappa. \qedhere
	\end{equation*}
\end{proof}

\begin{proof}[Proof of \Cref{thm:err without CIP}]
	For a coarse nodal basis function $\varphi_z \in V_H$, let $\varphi_z^{\rm ms, \gamma}$ and $\varphi_z^{\rm ms, 0}$ be the CIP and non-CIP finite element approximations, respectively, of the continuous saddle-point solution $\varphi_z^{\rm ms}$. Let $\II_h \varphi_z^{\rm ms} \in V_h$ denote the nodal interpolant. Since the continuous Lagrange multiplier $\mathcal{S}^{-1}\varphi_z$ belongs to $V_H$, no approximation error arises from the multiplier space. The stability of the corresponding mixed problems, the standard interpolation estimate, and $\kappa h \lesssim 1$ give 
	\begin{equation*}
		 \norm{\varphi_z^{\rm ms, \gamma} - \varphi_z^{\rm ms, 0}}_\gamma \leq \norm{\varphi_z^{\rm ms, \gamma} - \II_h \varphi_z^{\rm ms}}_\gamma + \norm{\varphi_z^{\rm ms, 0} - \II_h \varphi_z^{\rm ms}}_\gamma \lesssim h \abs{\varphi_z^{\rm ms}}_{H^2}.
	\end{equation*}
	Using \Cref{lem:H2}, we obtain
	\begin{equation*}
		\norm{\IdC \varphi_z - \mathcal{C}_{\infty}^0 \varphi_z}_\gamma = \norm{\varphi_z^{\rm ms, \gamma} - \varphi_z^{\rm ms, 0}}_\gamma \lesssim h \abs{\varphi_z^{\rm ms}}_{H^2} \lesssim (1 + \frac1{\kappa H})\kappa h.
	\end{equation*}
	Next, \Cref{thm:corr-appr} applied to the non-CIP corrector yields
	\begin{equation*}
		\norm{(\mathcal{C}_{\infty}^0 - \mathcal{C}_L^0)\varphi_z}_\gamma \le C_{\rm glob}^{1/2} e^{-\delta L / 2\tau}. 
	\end{equation*}
	Therefore,
	\begin{equation*}
		\norm{(\IdC - \mathcal{C}_{\infty}^0)\varphi_z}_\gamma \le \norm{(\IdC - \mathcal{C}_{\infty}^0) \varphi_z}_\gamma + \norm{(\mathcal{C}_{\infty}^0 - \mathcal{C}_L^0)\varphi_z}_\gamma \lesssim (1 + \frac1{\kappa H})\kappa h + C_{\rm glob}^{1/2} e^{-\delta L / 2\tau}. 
	\end{equation*}
	By linearity, the stability of the coarse nodal basis, and the finite-overlap argument, the preceding estimate for a single nodal basis function extends to arbitrary $v_H \in V_H$. Collecting the resulting factors depending on $\kappa H$, the PML parameter $\sigma_0$, and the localization parameter $L$ into $C(\kappa H, \sigma_0, L)$, we obtain
	\begin{equation*}
		\norm{(\IdC - \mathcal{C}_{\infty}^0)v_H}_\gamma \le C(\kappa H, \sigma_0, L) (\kappa h + e^{-\delta L / 2\tau}) \norm{v_H}_\gamma.
	\end{equation*}
	The same argument applies to the adjoint correctors and yields
	\begin{equation*}
		\norm{(\aIdC - \mathcal{C}_{\infty}^{0,\star})v_H}_\gamma \le C(\kappa H, \sigma_0, L) (\kappa h + e^{-\delta L / 2\tau}) \norm{v_H}_\gamma.
	\end{equation*}
	Set $\varepsilon_L := C(\kappa H, \sigma_0, L)(\kappa h + e^{-\delta L / 2\tau})$. For $v_{\infty}^{\rm ms} = (\II - \IdC) v_H$, we have $Q_H v_{\infty}^{\rm ms} = v_H$ and $(\II - \mathcal{Q}_L^{\rm ms}) v_{\infty}^{\rm ms} = (\mathcal{C}_{\infty}^0 - \mathcal{C}_{\infty}^\gamma) v_H$. Absorbing the stability constant of $Q_H$ into $C(\kappa H, \sigma_0, L)$ implies $\norm{(\II - \mathcal{Q}_L^{\rm ms})v_{\infty}^{\rm ms}} _\gamma \le \varepsilon_L \norm{v_{\infty}^{\rm ms}}_\gamma$ and the analogous estimate for the adjoint transfer operator. Thus \Cref{asp:abs-decay} is satisfied. The remaining assumptions of the abstract framework follow directly from the preservation of the coarse degrees of freedom and the stability of the ideal multiscale spaces. Therefore, under the smallness condition stated in the theorem, \Cref{thm:abs-loc-infsup} yields the stability of the improved trail-test pair, whereas \Cref{thm:err} gives
	\begin{equation*}
		\norm{u_{H, \infty, {\rm ms}}^\gamma - \pmsIII{u}{0}}_\gamma \le C(\kappa H, \sigma_0, L) (\kappa h + e^{-\delta L / 2\tau}) \norm{f}.
	\end{equation*}
	Combining this estimate with the standard localization bound completes the proof.
\end{proof}

\section{Numerical Experiments}\label{sec:numeric}
This section investigates the numerical stability and approximation properties of the improved CIP-LOD method proposed in \Cref{sec:alg-opt}. We consider a tow-dimensional Helmholtz radiation problem, for which the outgoing free-space solution is available in closed form. This allows us first to assess the combined PML truncation and fine-grid discretization error against the exact solution, and subsequently to isolate the coarse-scale LOD error by using the corresponding fine-grid solution as the reference.

Let $R_0 = 0.4$, $d_{\rm pml} = 0.1$, and $R = R_0 + d_{\rm pml} = 0.5$. The computational domain is $\Omega = B_R = (-R, R)^2$, whereas the physical domain is $B_{R_0} = (-R_0, R_0)^2$. All errors reported below are evaluated over $B_{R_0}$, where the PML coefficients reduce to the physical coefficients. The source is radially symmetric and compactly supported in the disk of radius $R_0$:
\begin{equation*}
	f^{\rm phys}(x) = \begin{cases}
		f_0, \quad &\abs{x} < R_0, \\
		0, \quad &\abs{x} \ge R_0,
	\end{cases}
\end{equation*}
Writing $r = \abs{x}$, the corresponding outgoing free-space solution is
\begin{equation*}
	u^{\rm ext}(r) = \begin{cases}
		\frac{\pi R_0 f_0}{2\kappa} \ii {\rm H}_1^{(1)}(\kappa R_0) {\rm J}_0(\kappa r) - \frac{f_0}{\kappa^2}, \quad & r < R_0, \\
		\frac{\pi R_0 f_0}{2\kappa} \ii {\rm J}_1(\kappa R_0) {\rm H}_0^{(1)}(\kappa r), \quad & r \ge R_0.
	\end{cases}
\end{equation*}
Here $\rm J_\nu$ and $\rm H_{\nu}^{(1)}$ denote the Bessel function of the first kind and the Hankel function of the first kind, respectively. 

A Cartesian PML is used to truncate the exterior domain. We prescribe the target reflection parameter $\sigma_{\max} = - \frac{3 \ln(R_{\rm pml})}{d_{\rm pml}}$, and define, for $j = 1, 2$,
\begin{equation*}
	\alpha_j (x_j) = \begin{cases}
		1 + \ii \frac{\sigma_{\max}}{\kappa} \big( \frac{\abs{x_j} - R_0}{d_{\rm pml}} \big)^2, \quad & \abs{x_j} > R_0, \\
		1, \quad & \abs{x_j} \le R_0.
	\end{cases}
\end{equation*}
The PML coefficients are $A(x) = \left( \begin{smallmatrix} \alpha_2(x_2)/\alpha_1(x_1) & 0\\ 0 & \alpha_1(x_1)/\alpha_2(x_2) \end{smallmatrix} \right)$ and $J(x) = \alpha_1(x_1)\alpha_2(x_2)$. A homogeneous Dirichlet condition is imposed on the outer boundary $\Gamma_{\rm out} = \partial B_R$. All discretization employ continuous piecewise affine finite elements. The coarse mesh $\TT_H$ is a uniform triangulation consisting of isosceles right triangles, and the fine mesh $\TT_h$ is obtained by uniform refinement of $\TT_H$. 

Let $u := u^{\rm ext}|_{B_{R_0}}$ denote the exact outgoing solution restricted to the physical domain $B_{R_0}$. For a numerical approximation $v$, we measure its relative error with respect to the exact solution by 
\begin{equation*}
	E(v) = \frac{\big( \kappa^2 \norm{u - v}_{B_{R_0}}^2 + \norm{\nabla(u - v)}_{B_{R_0}}^2 \big)^{1/2}}{\big( \kappa^2 \norm{u}_{B_{R_0}}^2 + \norm{\nabla u}_{B_{R_0}}^2 \big)^{1/2}}.
\end{equation*}
When a stabilized fine-grid solution $u_h^\gamma$ is used as the reference, we define
\begin{equation*}
	E_h(v) = \frac{\big( \kappa^2 \norm{u_h^\gamma - v}_{B_{R_0}}^2 + \norm{\nabla(u_h^\gamma - v)}_{B_{R_0}}^2 \big)^{1/2}}{\big( \kappa^2 \norm{u_h^\gamma}_{B_{R_0}}^2 + \norm{\nabla u_h^\gamma}_{B_{R_0}}^2 \big)^{1/2}}.
\end{equation*}
In comparisons between the stabilized and unstabilized formulations, $u_h$ always denotes the fine-grid solution obtained using the same penalty parameter as the corresponding multiscale approximation.

All multiscale approximations considered below incorporated a localized source correction. For each $T\in \faces_2(\TT_H)$, let $q_{T,L}^\gamma \in W_h(N^L(T))$ be defined by 
\begin{align}\label{sourceeee}
	a_h^\gamma(q_{T,L}^\gamma, w_h) = (f, w_h)_T \quad \forall w_h \in W_h(N^L(T)),
\end{align}
and set $q_L^\gamma = \sum_{T \in \faces_2(\TT_H)} q_{T, L}^\gamma$. Accordingly, we compute $\pmsI{\tilde  u} \in \pmsI{V}$ from 
\begin{align}
	a_h^\gamma(\pmsI{\tilde u}, \apmsI{v}) = (f, \apmsI{v}) - a_h^\gamma(q_L^\gamma, \apmsI{v}) \quad \forall \apmsI{v} \in \apmsI{V}.
\end{align}
The final source-corrected approximation is $\pmsIII{u}{\gamma, {\rm src}} = \pmsI{\tilde u} + q_L^\gamma$.

\subsection{Approximation properties and pollution reduction}
We now investigate the approximation properties of the CIP-LOD method. The penalty parameters used in this subsection depend on the purpose of the individual experiment. In the first experiment, which studies convergence with respect to the coarse mesh size, we employ the uniform purely imaginary penalty $\gamma_e = -0.01\ii$. All subsequent experiments employ the dispersion-optimized real penalty parameters. More precisely, $\gamma_e = -0.084$ on the horizontal and vertical fine-grid edges, which form the legs of the isosceles right triangles, and $\gamma_e = -0.041$ on the diagonal edges corresponding to their hypotenuses.

We distinguish the following four approximations. The standard coarse finite element solution is denoted by $u_H \in V_H$. The unstabilized LOD approximation $u_{H, L, {\rm ms}}$ is obtained by omitting the CIP penalty from both the local corrector problems and the global Petrov-Galerkin system. The fully stabilized CIP-LOD approximation $\pmsI{u}$ uses the form $a_h^\gamma$ in both the local and global problems. Finally, $\pmsIII{\hat{u}}{0}$ denotes the improved construction introduced in \Cref{sec:alg-opt}: its local correctors are computed with $a_h^0$, whereas its global Petrov-Galerkin system is assembled with $a_h^\gamma$. Similarly, $\hat{u}_{H,L,ms}^{0,src}$ the improved source-term-corrected construction: as with the element correctors, no penalty is added in the computation of Example \ref{sourceeee}. Then, $u_{H,L,ms}^{src}$ denotes LOD-FEM with source-corrected approximation. The individual experiments display the subset of these approximations relevant to the property under consideration.

\begin{example}[Convergence with respect to the coarse mesh size]
We first study convergence with respect to the coarse mesh size. We fix $\kappa = 32$ and $h = 2^{-10}$, and consider the sequence of coarse meshes $H = 2^{-j}, j = 4,5,...,8$. For each coarse mesh, the number of localization layers is chosen as $L = j-1$, so that $L$ increases from $3$ to $7$ as the coarse mesh is refined. \Cref{fig:order} shows that the source-corrected multiscale approximation converge steadily toward the stabilized fine-grid solution. Over the considered refinement range, their relative errors exhibit an observed convergence rate of approximately $E_h \sim H^{1.5}$. This enhanced rate is consistent with the source-regularity estimate established above. The present source is piecewise constant. More precisely, $f = 1$ inside the central disk and $0$ outside. Consequently, $f \in H^s(\Omega)$ for every $s < 1/2$. The multiscale correction improves the accuracy relative to the coarse finite element method, and the source correction yields an additional reduction in the error. Moreover, the results obtained with $a_h^\gamma$ and $a_h^0$ in the local corrector problems are nearly indistinguishable. Thus, omitting the CIP jump term from the local problems does not lead to a perceptible loss of accuracy in this experiment. This observation supports the improved construction of \Cref{sec:alg-opt}.

\begin{figure}[htbp]
	\centering
	\includegraphics[width=0.8\textwidth]{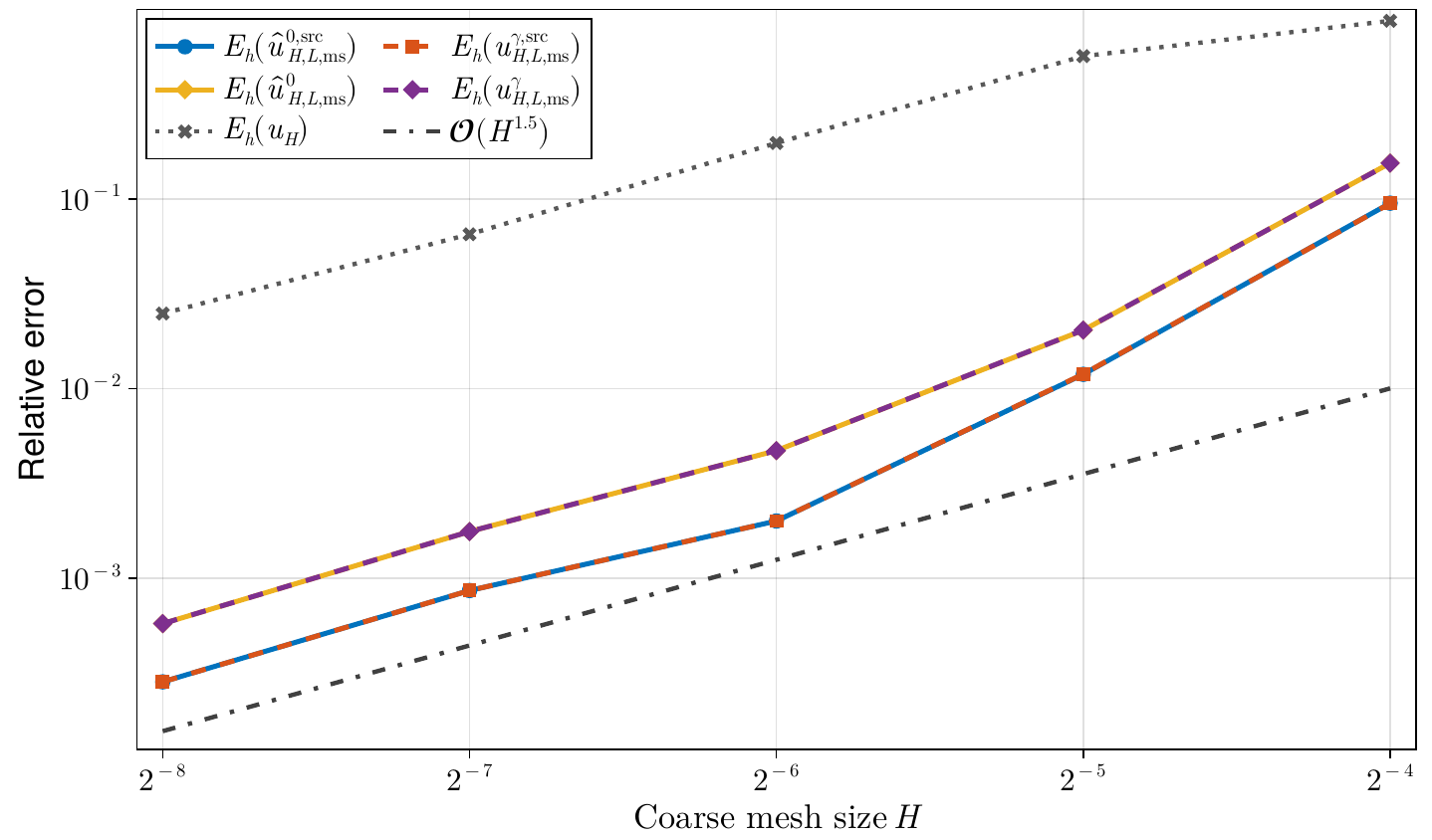}
	\centering
	\caption{Relative errors with respect to the stabilized fine-grid solution versus the coarse mesh size $H$.}
	\label{fig:order}
\end{figure}
\end{example}

\begin{example}[Pollution behavior at fixed resolution]
We next investigate the dependence of the multiscale error on the wave number under a fixed resolution condition. We impose $\kappa H = 2$ and $\kappa h = 1/4$, and vary $\kappa = 2^j, j = 5,6,...,10$. The oversampling parameter is chosen as $L = \lfloor \ln \kappa + 1/2 \rfloor - 1$, where $\lfloor x \rfloor$ denotes the greatest integer not exceeding $x$. \Cref{fig:pollution} shows that the standard coarse finite element solution has a relative error greater than one throughout the tested wave-number range. This behavior is expected because $\kappa H = 2$ lies outside the coarse-grid resolution regime. The multiscale correction substantially improves the coarse finite element solution, but the unstabilized LOD approximation still exhibits wave-number-dependent error growth. Its relative error increases progressively with $\kappa$ and reaches approximately one at $\kappa = 1024$. By contrast, both CIP-LOD variants remain stable throughout the tested range. Their relative errors stay close to $5\%$. Moreover, at large wave numbers, the errors obtained with and without the CIP penalty in the local corrector problems are almost indistinguishable. These results support the improved construction of \Cref{sec:alg-opt}.

\begin{figure}[htbp]
	\centering
	\includegraphics[width=0.8\textwidth]{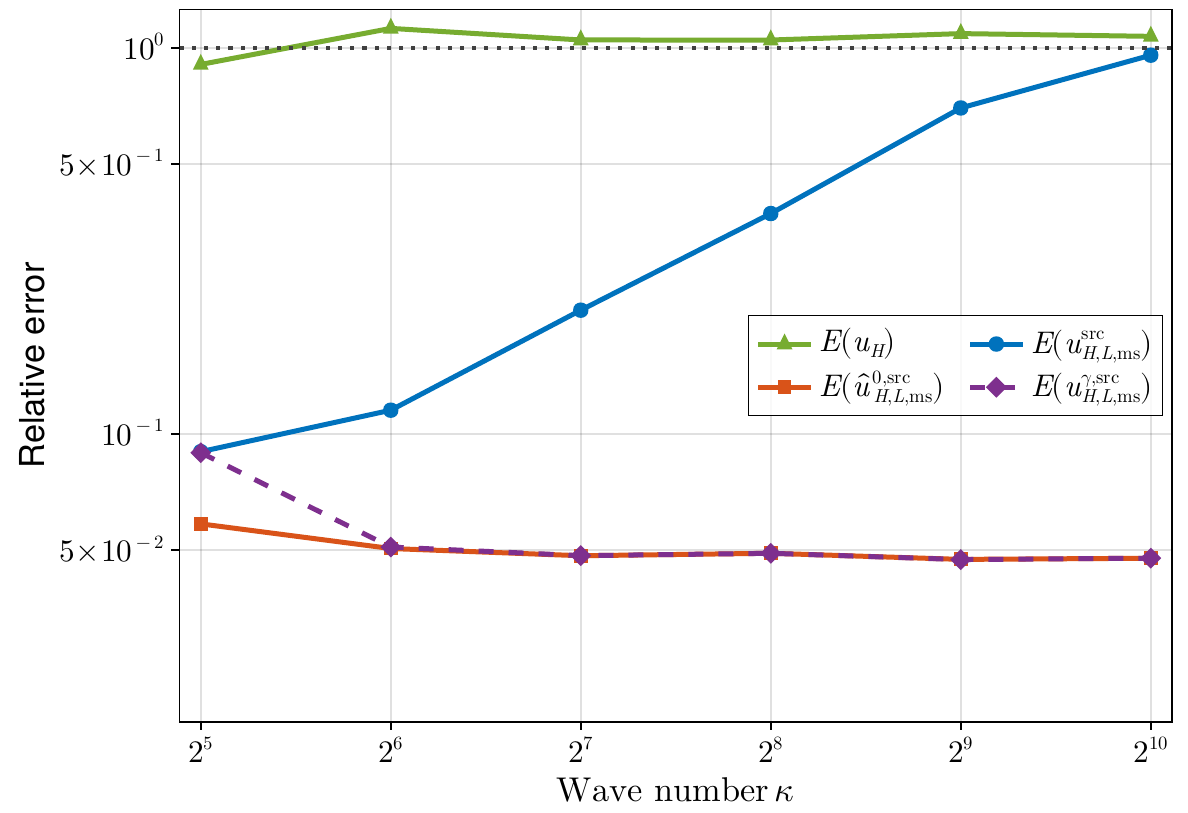}
	\caption{Relative errors versus the wave number $\kappa$ under the fixed-resolution condition $\kappa H = 2$ and $\kappa h = 1/4$.}
	\label{fig:pollution}
\end{figure}
\end{example}

\begin{example}[Influence of the number of the localization layers]
Finally, we examine the influence of the oversampling parameter in the high-frequency regime. We fix $\kappa=1024$, $\kappa H = 2$, and $\kappa h = 1/4$. Then we vary the number of localization layers from $L = 1$ to $L = 6(=\lfloor \ln \kappa + 1/2 \rfloor -1)$. \Cref{fig:patch} shows that enlarging the local patches does not improve the unstabilized LOD approximation. Even when the number of localization layers is increased to $L=6$, its relative error remains close to one and exhibits no meaningful decay. This indicates that the dominant error is not cause by localization. In contrast, the errors of both CIP-LOD variants decrease rapidly as $L$ increases. At $L=5$, both method attain a relative error of approximately $5\%$. Furthermore, the two CIP-LOD approximations exhibit essentially the same decay with respect to $L$ and reach the same final accuracy. This confirms that decisive stabilization is provided by the optimized penalty in the global Petrov-Galerkin system, whereas including the CIP jump term in the local corrector problems offers no visible additional benefit.

\begin{figure}[htbp]
	\centering
	\includegraphics[width=0.8\textwidth]{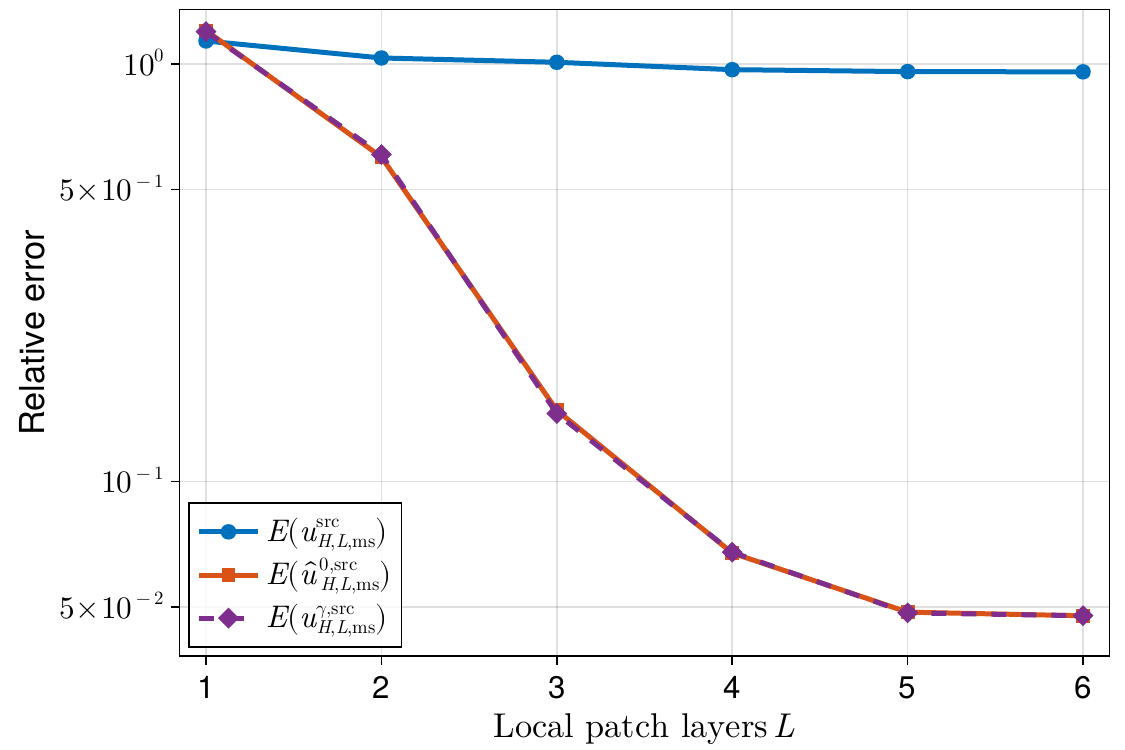}
	\caption{Relative errors versus the oversampling parameter $L$ for $\kappa = 1024, \kappa H = 2$, and $\kappa h = 1/4$.}
	\label{fig:patch}
\end{figure}
\end{example}

\subsection{Exploratory test beyond the theoretical coarse-resolution regime}\label{sec:exp-stab}
We next explore the behavior of the CIP-LOD method in a coarse-grid regime, which is not covered by the present analysis. In particular, we allow $\kappa H$ to take values substantially larger than those required by the theoretical results above. The purpose of this experiment is purely numerical: we aim to assess how far the stabilized multiscale method can remain reliable under severe coarse-grid under-resolution. A rigorous explanation of this behavior, in particular the stability mechanism for large \(\kappa H\), is an interesting direction for future work.

We fix $\kappa = 960$, $\kappa h =1$, and $\gamma_e = -0.1\ii$, and vary the coarse mesh size such that $\kappa H = 2,...,12$. To investigate the influence of localization, the experiment is repeated for $L = 3,...,7$. Both the standard LOD approximation $u_{H, L, {\rm ms}}$ and the CIP-LOD $\pmsI{u}$ are compared with the same stabilized fine-grid solution. 

\begin{figure}[htbp]
	\centering
	\includegraphics[width=0.8\textwidth]{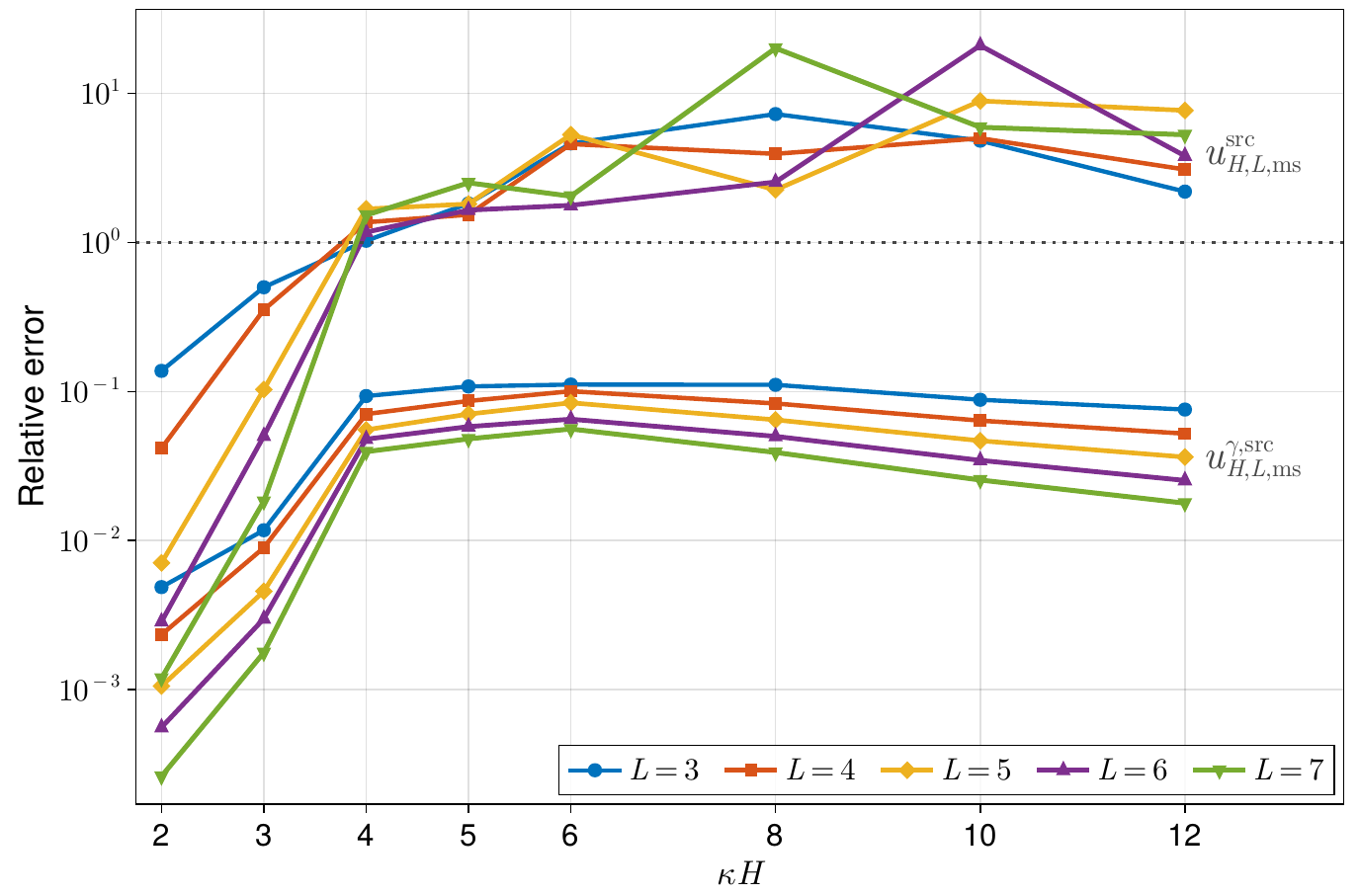}
	\centering
	\caption{Relative error of the CIP-FEM and CIP-LOD solutions with respect to their corresponding stabilized fine-grid solution.}
	\label{fig:kH}
\end{figure}

\Cref{fig:kH} reveals a pronounced difference between the two multiscale formulations. For the standard LOD method, the relative error exceeds one when $\kappa H \ge 4$ and for every tested value of $L$. Increasing the number of localization layers does not produce a meaningful reduction in the error. By contrast, the CIP-LOD approximation remains  numerically robust through the tested range. Its relative error stays close to $10\%$ even for $\kappa H = 12$. This behavior suggests that the stabilizing effect of the CIP term persists far beyond the parameter regime covered by the current theory. However, we emphasize that this observation is empirical; establishing stability and localization results for such strongly under-resolved coarse meshes remains an open analytical question and is left for future work.

\section{Conclusion}\label{sec:con}
This paper develops a CIP-LOD framework for high-frequency Helmholtz radiation problems. The method combines PML truncation for the exterior problem, CIP stabilization for the fine-scale discretization, and localized orthogonal decomposition for the construction of multiscale Petrov--Galerkin spaces on coarse meshes. Starting from the CIP-PML reference discretization, we show that, under a coarse-scale resolution condition of the form $\kappa H \lesssim 1$, the sesquilinear form is coercive on the fine-scale kernel associated with the quasi-interpolation operator. This yields the well-posedness of the corrector problems and exponential decay of the corresponding element correctors. Consequently, the localized CIP correctors converge exponentially to their ideal counterparts as the oversampling parameter increases. We further discuss the enhanced stability provided by the CIP stabilization and its influence on the stability and localization of the corrector problems.

To analyze localized spaces constructed with a perturbed corrector operator, we introduce an abstract stability-transfer framework. The framework shows that a localized trial-test pair inherits the stability of the ideal Petrov--Galerkin formulation if it preserves the same coarse degrees of freedom and approximates the ideal multiscale spaces sufficiently accurately. This result is then applied to the improved construction proposed in this paper, in which the CIP penalty is omitted from the local corrector problems while retained in the global formulation. The resulting analysis shows that the perturbation introduced by omitting the local CIP term is controlled by the fine-grid discretization and localization errors.

The numerical experiments confirm the approximation properties of the proposed method. The CIP-LOD approximation substantially improves the coarse finite element solution and effectively suppresses the pollution error. Moreover, the fully stabilized and improved CIP-LOD constructions exhibit nearly indistinguishable accuracy in the tested regimes. It supports the use of the cheaper local corrector problems without CIP penalties. The experiments also reveal that the stabilized multiscale method remains numerically robust for values of $\kappa H$ far beyond the regime covered by the present analysis. This observation is purely numerical and suggests that the theoretical coarse-scale resolution condition may be conservative for the stabilized formulation. A rigorous explanation of this behavior is left for future work.

Although the present analysis is restricted to a two-dimensional Helmholtz radiation problem with Cartesian PML truncation, the stability-transfer principle itself is formulated independently of this particular setting. When the corresponding stability, interpolation, regularity, and localization estimates are available for the underlying problem, the same construction can be carried over to three-dimensional problems, formulations with impedance boundary conditions, and heterogeneous media. Hence, the CIP-LOD strategy and the improved corrector construction are not intrinsically tied to the particular two-dimensional model considered here.

\section*{Acknowledgements}
This work was supported by National Key Research and Development Project of China (Grant No. 2024YFA1012600). Huangxin Chen was also supported by National Natural Science Foundation of China (Grant No. 12471345, 12526102) and Fujian Provincial Natural Science Foundation of China (Grant No. 2026J011002).


\end{document}